\input ./style/arxiv-general.cfg
\documentclass[aos,MSNbibl,seceqn,nameyear,dvips]{arximspdf}
\makeatletter
   \@ifpackageloaded{graphicx}{}{\usepackage{graphicx}}
\makeatother
\usepackage{accents}
\usepackage{upgreek}

\doi{10.1214/15-AOS1356}% Updated by VTEXPTS2LaTeX.exe, 16.07.2015
\volume{43}
\issue{6}
\pubyear{2015}
\firstpage{2676}
\lastpage{2705}
\docsubty{FLA}

\makeatletter
\newcommand{\rrvert}{\vert}

\newcommand{\llvert}{\vert}
\newcommand{\T}{ \mathrm{T}}
\newcommand{\bbeta}{\bolds{\beta}}
\newcommand{\btheta}{\bolds{\theta}}
\newcommand{\bphi}{\bolds{\phi}}
\newcommand{\0}{\mathbf{0}}
\newcommand{\ba}{\mathbf{a}}
\newcommand{\bb}{\mathbf{b}}
\def\T{\mathrm{T}}
\def\balpha{\bolds{\alpha}}
\def\bC{\mathbf{C}}
\newcommand{\GIC}{\mathrm{GIC}}
\newcommand{\NO}{\mathrm{NO}}
\newcommand{\SO}{\mathrm{SO}}
\newcommand{\calO}{\mathcal{O}}
\newcommand{\eeta}{\bolds{\eta}}
\newcommand{\bT}{\mathbf{T}}
\newcommand{\bH}{\mathbf{H}}
\newcommand{\calL}{\mathcal{L}}
\newcommand{\calA}{\mathcal{A}}
\newcommand{\calB}{\mathcal{B}}
\newtheorem{theorem}{Theorem}[section]
\newproclaim{rem}{Remark}[section]
\newtheorem{prop}{Proposition}[section]
\newproclaim{exa}{Example}[section]
\newproclaim{ass}{Assumption}
\makeatother

\begin{document}
\begin{frontmatter}

\title{Model selection and structure specification in ultra-high
dimensional generalised semi-varying coefficient models\thanksref{T1}}
\runtitle{Model selection in GSVCM}
\thankstext{T1}{Supported by the
Singapore National Research Foundation under its Cooperative Basic Research
Grant and administered by the Singapore Ministry of Health National Medical
Research Council Grant No. NMRC/CBRG/0014/2012 and the National Science
Foundation of China Grant No.~11271242.}

\begin{aug}
% Corresponding author: Wenyang Zhang - wenyang.zhang@york.ac.uk% Updated by VTEXPTS2LaTeX.exe, 17.07.2015 08:02
%Updated by VTEXPTS2LaTeX.exe, 16.07.2015 10:52
\author[A]{\fnms{Degui}~\snm{Li}\ead[label=e1]{degui.li@york.ac.uk}},
\author[A]{\fnms{Yuan}~\snm{Ke}\ead[label=e2]{yk612@york.ac.uk}}
\and
\author[A]{\fnms{Wenyang}~\snm{Zhang}\corref{}\ead[label=e3]{wenyang.zhang@york.ac.uk}}
\runauthor{D. Li, Y. Ke and W. Zhang}
%\runauthor{}
\affiliation{University of York}
%\dedicated{}
\address[A]{Department of Mathematics\\
University of York\\
Heslington\\
York YO10 5DD\\
United Kingtom\\
\printead{e1}\\
\phantom{E-mail:\ }\printead*{e2}\\
\phantom{E-mail:\ }\printead*{e3}}
\end{aug}

% HISTORY:
%
\received{\smonth{1} \syear{2015}}% Updated by VTEXPTS2LaTeX.exe,
%16.07.2015 10:52
%
\revised{\smonth{6} \syear{2015}}% Updated by VTEXPTS2LaTeX.exe,
%16.07.2015 10:52

% ABSTRACT
%
\begin{abstract}

In this paper, we study the model selection and structure specification
for the generalised semi-varying coefficient models (GSVCMs), where the
number of potential covariates is allowed to be larger than the sample
size. We first propose a penalised likelihood method with the LASSO
penalty function to obtain the preliminary estimates of the functional
coefficients. Then, using the quadratic approximation for the local
log-likelihood function and the adaptive group LASSO penalty (or the
local linear approximation of the group SCAD penalty) with the help of
the preliminary estimation of the functional coefficients, we introduce
a novel penalised weighted least squares procedure to select the
significant covariates and identify the constant coefficients among the
coefficients of the selected covariates, which could thus specify the
semiparametric modelling structure. The developed model selection and
structure specification approach not only inherits many nice
statistical properties from the local maximum likelihood estimation and
nonconcave penalised likelihood method, but also computationally
attractive thanks to the computational algorithm that is proposed to
implement our method. Under some mild conditions, we establish the
asymptotic properties for the proposed model selection and estimation
procedure such as the sparsity and oracle property. We also conduct
simulation studies to examine the finite sample performance of the
proposed method, and finally apply the method to analyse a real data
set, which leads to some interesting findings.
\end{abstract}

% KEYWORDS
%
%\begin{keyword}[class=AMS]
%\kwd[Primary ]{}
%\kwd{}
%\kwd[; secondary ]{}
%\end{keyword}
%\begin{keyword}
%\kwd{}
%\end{keyword}
\end{frontmatter}

%%%%%%%%%%%%%%%%%%%%%%

%s1 #&#
\section{Introduction}
\label{intro}

In recent years, model selection has become an important and
fundamental issue in data analysis as high-dimensional data are
commonly encountered in various applied fields such as epidemiology,
genetics and finance. It is well known that the traditional model
selection procedures such as the stepwise regression and the best
subset variable selection can be extremely computationally intensive in
the analysis of the high-dimensional data. To address this
computational challenge, various penalised likelihood/least-square
methods have been well studied and become a promising alternative. With
an appropriate penalty function, the penalised method would
automatically shrink the small coefficients to zero and remove the
associated variables from the model, hence serve the purpose of model
selection. Some commonly-used penalty functions include the LASSO
[\citet{Tib96}], SCAD [\citet{FanLi01}], group LASSO [\citet{YuaLin06}], adaptive LASSO [\citet{Zou06}] and MCP [\citet{Zha10}], and the
algorithms to implement the penalised likelihood/least squares methods
have also been developed in the literature [cf., \citet{Efretal04,HunLi05,ZouLi08}]. In high-dimensional data
analysis, it is often the case that the number of potential covariates
grows over sample size or even diverges with certain exponential rate.
In the context of parametric models, there has been some literature
addressing this problem; see, for example, \citet{HuaXie07}, \citet{FanLv08}, Huang, Horowitz and Ma (\citeyear{HuaHorMa08}), \citet{ZhaHua08},
Fan, Samworth and Wu (\citeyear{FanSamWu09}), \citet{ZouZha09}, \citet{FanSon10} and
B\"
{u}hlmann and van de Geer (\citeyear{Buhvan11}).\looseness=1

However, the pre-supposed parametric linear relationships and models,
although easy to implement, are often too restricted and unrealistic in
practical application. They often lead to model misspecification, which
would result in inconsistent estimates and incorrect conclusions being
drawn from the data analysed. In this paper, we relax this linear
restriction and use functional coefficients to describe the
relationship between response and covariates. The varying coefficient
models, as an important and useful generalisation of the linear models,
have played a very important role in the analysis of the complex data
and experienced deep and exciting developments; see, for example, Fan
and Zhang (\citeyear{FanZha99,FanZha00}), Cheng, Zhang and Chen (\citeyear{CheZhaChe09}), \citet{WanXia09},
Wang, Kai and Li (\citeyear{WanKaiLi09}), Zhang, Fan and Sun (\citeyear{ZhaFanSun09}), Kai, Li and Zou
(\citeyear{KaiLiZou11}) and \citet{LiZha11}. Suppose we have a response variable
$y$, an index variable $U$, and potential covariates $X=(x_1, \ldots,
x_{d_n})^{\T}$, where the dimension $d_n$ depends on sample size $n$
and $d_n \rightarrow\infty$ when $n \rightarrow\infty$. Define the
conditional expectation of $y$ for given $(U, X)$ by
\[
m(U, X) = \mathsf{ E}(y| U, X).
\]
We assume in this paper that the log conditional density function of
$y$ given $X$ and $U$ is
%
%e1.1 #&#
\begin{equation}\qquad
C_1(\bphi_1) \ell \bigl(m(U, X), y \bigr) +
C_2(y, \bphi_2) \qquad\mbox{with } g \bigl(m(U, X) \bigr) = \sum
_{j=1}^{d_n} a_j(U)
x_j, \label{eq1.1}
\end{equation}
where $g(\cdot)$, $\ell(\cdot, \cdot)$, $C_1(\cdot)$ and
$C_2(\cdot, \cdot)$ are known, the functional coefficients $a_1(\cdot), \ldots,
a_p(\cdot)$ are unknown to be estimated, $C_1(\bphi_1) > 0$, $\bphi_1$
and $\bphi_2$ are unknown nuisance parameters. When the response
variable is discrete, we define the density function as its probability
mass function. It is easy to see that model (\ref{eq1.1}) is a natural
extension of the generalised linear models by allowing the coefficients
varying with the index variable $U$. As some functional coefficients in
(\ref{eq1.1}) may be constant coefficients, we call (\ref{eq1.1}) as
generalised semi-varying coefficient models (GSVCMs).

The model selection in the varying coefficient models (which can be
seen as a
special case of the GSVCMs) has been extensively studied in existing
literature. For instance, Wang, Li and Huang (\citeyear{WanLiHua08}) and \citet{WanXia09} use the group
penalisation to select the significant variables in the varying coefficient
models when the number of potential covariates is fixed. More recently, for
the ultra-high dimensional varying coefficient models, Song, Yi and Zuo (\citeyear{SonYiZuo12}),
\citet{Cheetal14}, Fan, Ma and Dai (\citeyear{FanMaDai14}) and Liu, Li and Wu (\citeyear{LiuLiWu14})
combine the nonparametric
independence screening technique and the group penalised method to
choose the
significant covariates and estimate the functional coefficients for the varying
coefficient models. Fan, Feng and Song (\citeyear{FanFenSon11}) consider a nonparametric
independence screening in sparse ultra-high dimensional additive models.
Wei, Huang and Li (\citeyear{WeiHuaLi11}) consider the penalised variable selection by
using a
basis function approximation for the functional coefficients in the varying
coefficient models and allows that the number of covariates diverges
with the
sample size. \citet{Lia12} further generalises Wei, Huang and Li's (\citeyear{WeiHuaLi11})
methodology to the generalised varying coefficient models which are
similar to
our framework (\ref{eq1.1}). Unlike the existing literature, in this paper,
the model selection for the proposed GSVCMs has two aspects: (1) variable
selection; and (2) identification of the constant coefficients. As the
variable selection is equivalent to identifying the zero functional
coefficients, and the identification of the constant coefficients is
equivalent to identifying the functional coefficients with zero
derivative or
variation. Either of the two aspects would be related to the so-called
``all-in-all-out'' problem.

In this paper, we first propose a penalised likelihood method with the LASSO
penalty function to obtain the preliminary estimates of the functional
coefficients, which is proved to be uniformly consistent. The uniform
convergence rate for the preliminary penalised nonparametric estimation
results relies
on the number of nonzero functional coefficients and the tuning parameter
involved in the penalty term. Then we use the preliminary estimates of the
functional coefficients in the quadratic approximation for the local
log-likelihood function and the construction of the adaptive group LASSO
penalty or the local linear approximation of the group SCAD penalty,
and introduce a novel penalised weighted least squares procedure to
simultaneously select the significant covariates and
identify the constant coefficients among the coefficients of the selected
covariates. Hence, the semi-varying coefficient modelling structure can be
specified. The developed model selection and structure specification approach
inherits many nice statistical properties from both the local maximum
likelihood estimation and nonconcave penalised likelihood method.
Under some
regularity conditions, we establish the asymptotic properties for the proposed
model selection and estimation procedure such as the sparsity and oracle
property. In order to implement our method in practice, we develop a novel
computational algorithm to do the maximisation involved in the estimation
procedure when the SCAD or LASSO penalty is used. The SCAD has many
advantages, and is
widely used as a penalty function in the shrinkage methods. The commonly used
approach, to deal with the SCAD penalty in the implementation of
shrinkage method for the varying coefficient models, consists of two
steps: (1) approximate SCAD by an $L_1$
penalty locally by local linear approximation; (2) apply the quadratic
approximation to deal with the $L_1$ penalty. In this paper, we do not go
down that route. Making use of the structure of the group SCAD, we
propose a
different algorithm to implement our method. Our simulation results
show that both the adaptive group LASSO and the SCAD methods perform
reasonably well with the latter giving slightly better
performance, and the method developed in the present paper outperforms those
in \citet{WanXia09}, and \citet{Lia12}.

The rest of the paper is organised as follows. Section~\ref{sec2}
describes the penalised model selection and structure specification
procedure. Section~\ref{sec3} gives the asymptotic properties of the
proposed model selection and structure specification procedure. Section~\ref{sec4}
provides a computational algorithm to implement the
developed method and discusses how to determine the tuning parameters.
Section~\ref{sec5} compares the finite sample performance of the
developed model selection with those proposed in the existing
literature through some simulation studies. In Section~\ref{sec6}, we
apply the GSVCMs together with the proposed model selection, structure
specification and estimation procedure to analyse an environmental data
set from Hong Kong, and explore how some pollutants and other
environmental factors affect the number of daily total hospital
admissions for circulationary and respirationary problems in Hong Kong.
The regularity conditions for the asymptotic theory are given in
Appendix \ref{app.A}. The proofs of the main theoretical results and
some auxiliary results are provided in Appendices B and C of a
supplemental document [Li, Ke and Zhang (\citeyear{supp})].

%%%%%%%%%%%%%%%%%%%%%%

%s2 #&#
\section{Model selection and structure specification method}
\label{sec2}

For any function $f(\cdot)$, throughout this paper, we use $\dot
{f}(\cdot)$ to denote its first-order derivative, and $\ddot{f}(\cdot)$
its second-order derivative. For any vector ${\mathbf{u}}$, we define
$\|{\mathbf{u}}\|^2={\mathbf{u}}^{\T}{\mathbf{u}}$. As in the
generalised linear models, our main interest lies in the conditional
mean of the response variable for given covariates, and $C_1(\bphi_1)$
and $C_2(y, \bphi_2)$ in model (\ref{eq1.1}) have little to do with the
mean part. In order to make the presentation simpler, without loss of
generality, we assume the log conditional density function of $y$ given
$X$ and $U$ is
%
%e2.1 #&#
\begin{equation}
\ell \bigl(m(U, X), y \bigr)\qquad \mbox{with } g \bigl(m(U, X) \bigr) = \sum
_{j=1}^{d_n} a_j(U) x_j
\label{eq2.1}
\end{equation}
and further assume the support of the index variable $U$ is $[0, 1]$
throughout this paper. Our model selection and structure specification
procedure can be summarised as follows: (i) use the penalised local
maximum likelihood method with the LASSO penalty to obtain the
preliminary estimation of the functional coefficients (see Section~\ref{sec2.1}); (ii) consider the quadratic approximation of the
log-likelihood function by using the preliminary functional
coefficients estimates and the approximated log-likelihood function is
essentially an $L_2$ objective function (see Section~\ref{sec2.2});
(iii) conduct the variable selection and structure specification by
using a penalised weighted least squares method with two types of
weighted group LASSO penalty functions: the adaptive group LASSO and
the local linear approximation of the group SCAD where the preliminary
functional coefficients estimates are also used (see Section~\ref{sec2.3}); (iv) finally estimate the constant coefficients in the
GSVCMs (see Section~\ref{sec2.4}). The model selection procedure
proposed in this paper can be seen, in some sense, as a generalisation
of Fan, Ma and Dai's (\citeyear{FanMaDai14}) folded concave penalised estimation for
ultra-high dimensional parametric regression models.

%s2.1 #&#
\subsection{Preliminary estimation of the functional
coefficients}\label{sec2.1}

Suppose we have a sample $(U_i, X_i, y_i)$, $i=1, \ldots, n$, from
model (\ref{eq2.1}), where $X_i = (x_{i1},\break \ldots, x_{id_n})^{\T}$. For
each given $k$, $k=1, \ldots, n$, by Taylor's expansion of $a_j(\cdot
)$, $j=1,\ldots,d_n$, we have
\[
a_j(U_i) \approx a_j(U_k) +
\dot{a}_j(U_k) (U_i - U_k),
\]
when $U_i$, $i=1, \ldots, n$, are in a small neighbourhood of $U_k$.
This local linear approximation leads to the construction of the
following local log-likelihood function to estimate $a_j(U_k)$ and
$\dot
{a}_j(U_k)$, $j=1, \ldots, d_n$:
%
%e2.2 #&#
\begin{eqnarray}
\label{eq2.2}&& {\mathcal{L}}_{nk}(\ba_k,
\bb_k)
\nonumber
\\[-8pt]
\\[-8pt]
\nonumber
&&\qquad=\frac{1}{n}\sum_{i=1}^n
\ell \Biggl(g^{-1} \Biggl\{ \sum_{j=1}^{d_n}
\bigl[\alpha_{jk} + \beta_{jk} (U_i -
U_k) \bigr] x_{ij} \Biggr\}, y_i \Biggr)
K_h(U_i - U_k),
\end{eqnarray}
where $K(\cdot)$ is a kernel function, $h$ is a bandwidth, $K_h(\cdot
)=\frac{1}{h}K(\cdot/h)$,
\[
\ba_k = (\alpha_{1k}, \ldots, \alpha_{d_nk}
)^{\T},\qquad \bb_k = (\beta_{1k}, \ldots,
\beta_{d_nk} )^{\T}.
\]
When the dimension of the covariates is fixed, we may obtain the
solution which maximises the local log-likelihood function ${\mathcal
{L}}_{nk}(\cdot,\cdot)$ defined in (\ref{eq2.2}) and show that the
resulting nonparametric estimators are consistent [cf., Cai,
Fan and Li \citeyear{CaiFanLi00}; \citet{ZhaPen10}]. However, for the case of the
ultra-high dimensional GSVCMs, it would be difficult to obtain
satisfactory estimation by maximising ${\mathcal{L}}_{nk}(\cdot,\cdot)$
as the number of unknown nonparametric components involved exceeds the
number of observations. In order to address this issue, we next
introduce a penalised local log-likelihood method by adding an
appropriate penalty function to the above local log-likelihood function.

Without loss of generality, we assume that there exist $1\leq s_{n1}<
s_{n2}<d_n$ such that for $1\leq j\leq s_{n1}$, $a_j(\cdot)$ are the
functional coefficients with nonzero deviation; for $s_{n1}+1\leq
j\leq s_{n2}$, $a_j(\cdot)\equiv c_j$ are the constant coefficients;
for $s_{n2}+1\leq j\leq d_n$, $a_j(\cdot)\equiv0$. Moreover, we assume
that $s_{n2}$, although may diverge with the sample size, is much
smaller than the sample size $n$ and the dimension of the whole
covariates $d_n$. Hence, for any $k=1,\ldots,n$, the number of nonzero
elements in $\ba_{k0}= [a_1(U_k), \ldots, a_{d_n}(U_k)
]^{\T
}$ and $\bb_{k0}= [\dot{a}_1(U_k), \ldots, \dot
{a}_{d_n}(U_k)
]^{\T}$ is at most $s_{n1}+s_{n2}$. Define the penalised local
log-likelihood function with the LASSO penalty function as
%
%e2.3 #&#
\begin{equation}
\label{eq2.3} {\mathcal Q}_{nk}(\ba_k,
\bb_k)={\mathcal L}_{nk}(\ba_k,
\bb_k) -\lambda_1\sum_{j=1}^{d_n}|
\alpha_{jk}| -\lambda_2\sum_{j=1}^{d_n}h|
\beta_{jk}|,
\end{equation}
where $\lambda_1$ and $\lambda_2$ are two tuning parameters. We let
$(\widetilde{\ba}_k, \widetilde{\bb}_k)$ be the maximiser of
${\mathcal Q}_{nk}(\cdot,\cdot)$ and call it the preliminary estimator
of the functional coefficients $a_j(U_k)$'s and their derivatives $\dot
{a}_j(U_k)$s, $j=1,\ldots,d_n$.

We will show in Proposition~\ref{pr3.1} that the above preliminary estimator
obtained by the penalised local likelihood estimation with the LASSO
penalty is uniformly consistent. The preliminary estimates of the
functional coefficients will be used in the approximation of
log-likelihood function and the construction of weighted LASSO penalty
functions in our model selection and structure specification procedure;
see Sections~\ref{sec2.2} and \ref{sec2.3} below.

%s2.2 #&#
\subsection{Quadratic approximation of the log-likelihood
estimation}\label{sec2.2}

In the\break model selection and structure specification procedure for the
GSVCMs, we need to consider the following local log-likelihood function:
%
%e2.4 #&#
\begin{equation}
\label{eq2.4} {\mathcal{L}}_n({\mathcal{A}}, {\mathcal{B}})=\sum
_{k=1}^n{\mathcal L}_{nk}(
\ba_k, \bb_k),
\end{equation}
where ${\mathcal{A}}= (\ba_1^{\T}, \ldots, \ba_n^{\T
} )^{\T}$,
${\mathcal{B}}= (\bb_1^{\T}, \ldots, \bb_n^{\T}
)^{\T}$,
and ${\mathcal L}_{nk}(\ba_k, \bb_k)$ is defined in (\ref{eq2.2}). To
alleviate the computational burden for the optimisation of ${\mathcal
{L}}_n({\mathcal{A}}, {\mathcal{B}})$, we next introduce a simple
approximation.

Let
\[
\dot{\mathcal{L}}_{n}({\mathcal{A}}, {\mathcal{B}}) = \bigl[\dot{
\mathcal{L}}_{n1}^{\T}(\ba_1,\bb_1),
\ldots, \dot {\mathcal {L}}_{nn}^{\T}(\ba_n,
\bb_n) \bigr]^{\T}
\]
and
\[
\ddot{\mathcal{L}}_{n}({\mathcal{A}}, {\mathcal{B}}) ={\mathsf{diag}}
\bigl\{\ddot{\mathcal{L}}_{n1}(\ba_1,\bb_1),
\ldots, \ddot{\mathcal{L}}_{nn}(\ba_n,\bb_n)
\bigr\},
\]
where
\begin{eqnarray*}
\dot{\mathcal L}_{nk}(\ba_k,\bb_k)&= &
\frac{1}{n}\sum_{i=1}^nq_1
\Biggl\{ \sum_{j=1}^{d_n} \bigl[
\alpha_{jk} +\beta_{jk} (U_i - U_k)
\bigr] x_{ij}, y_i \Biggr\} \pmatrix{ X_i
\vspace*{2pt}\cr
\displaystyle\frac{U_i-U_k}{h}\cdot X_i }
\\
&&{} \times K_h(U_i-U_k),
\\
\ddot{\mathcal{L}}_{nk}(\ba_k,\bb_k, l) &= &
\frac{1}{n}\sum_{i=1}^n
q_2 \Biggl\{ \sum_{j=1}^{d_n}
\bigl[\alpha_{jk} +\beta_{jk} (U_i -
U_k) \bigr] x_{ij}, y_i \Biggr\} \biggl(
\frac{U_i-U_k}{h} \biggr)^{l}
\\
&&{} \times X_iX_i^{\T}K_h(U_i
- U_k),\qquad l=0,1,2,
\\
\ddot{\mathcal{L}}_{nk}(\ba_k,\bb_k)&= &
\left[\matrix{ \ddot{\mathcal{L}}_{nk}(
\ba_k,\bb_k, 0)&\ddot{\mathcal {L}}_{nk}(\ba
_k,\bb_k,1)
\vspace*{2pt}\cr
\ddot{\mathcal{L}}_{nk}(\ba_k,\bb_k, 1)&
\ddot{\mathcal {L}}_{nk}(\ba _k,\bb_k, 2)
}
 \right]
\end{eqnarray*}
and
\[
q_1(s,y)=\frac{\partial\ell [g^{-1}(s),y ]}{\partial s},\qquad q_2(s,y)=\frac{\partial^2 \ell [g^{-1}(s),y ]}{\partial s^2}.
\]
Denote $\widetilde{\mathcal{A}}_n= (\widetilde\ba_1^{\T},
\ldots,
\widetilde\ba_n^{\T} )^{\T}$ and
$\widetilde{\mathcal{B}}_n= (\widetilde\bb_1^{\T}, \ldots,
\widetilde
\bb_n^{\T}  )^{\T}$ where $(\widetilde{\ba}_k, \widetilde{\bb}_k)$
is the maximiser of the objective function ${\mathcal Q}_{nk}(\cdot
,\cdot)$ in (\ref{eq2.3}), and define
\[
{\mathcal V}_n({\mathcal{A}}, {\mathcal{B}})= \bigl(
\ba_1^{\T}, \bb _1^{\T
}, \ldots,
\ba_n^{\T}, \bb_n^{\T}
\bigr)^{\T},\qquad{\mathcal V}_n({\mathcal {A}}, h{\mathcal{B}})=
\bigl(\ba_1^{\T}, h\bb_1^{\T}, \ldots,
\ba _n^{\T}, h\bb_n^{\T}
\bigr)^{\T}.
\]

By Taylor's expansion of the log-likelihood function ${\mathcal
{L}}_n({\mathcal{A}}, {\mathcal{B}})$, we may obtain the following
quadratic approximation:
%
%e2.5 #&#
\begin{eqnarray}
\label{eq2.5} {\mathcal{L}}_n({\mathcal{A}}, {\mathcal{B}}) &
\approx& {\mathcal{L}}_{n}(\widetilde{\mathcal{A}}_n,
\widetilde {\mathcal{B}}_n)+ \bigl[{\mathcal V}_n({
\mathcal{A}}, h{\mathcal {B}})-{\mathcal V}_n(\widetilde{
\mathcal{A}}_n, h\widetilde{\mathcal {B}}_n)
\bigr]^{\T}\dot{\mathcal{L}}_{n}(\widetilde{
\mathcal{A}}_n, \widetilde{\mathcal{B}}_n)
\nonumber
\\
&&{}+\tfrac{1}{2} \bigl[{\mathcal V}_n({\mathcal{A}}, h{
\mathcal{B}})-{\mathcal V}_n(\widetilde{\mathcal{A}}_n, h
\widetilde{\mathcal{B}}_n) \bigr]^{\T
}
\nonumber
\\[-8pt]
\\[-8pt]
\nonumber
&&{}\times\ddot{
\mathcal{L}}_n(\widetilde{\mathcal{A}}_n, \widetilde {
\mathcal{B}}_n) \bigl[{\mathcal V}_n({\mathcal{A}}, h{
\mathcal{B}})-{\mathcal V}_n(\widetilde{\mathcal{A}}_n, h
\widetilde{\mathcal{B}}_n) \bigr]
\\
&\equiv&{\mathcal L}_n^\diamond({\mathcal{A}}, {\mathcal{B}}).\nonumber
\end{eqnarray}
It is easy to see that ${\mathcal{L}}_{n}^{\diamond}({\mathcal{A}},
{\mathcal{B}})$ is essentially an $L_2$ objective function. Hence, it
would be much easier to deal with ${\mathcal{L}}_{n}^{\diamond
}({\mathcal{A}}, {\mathcal{B}})$ in (\ref{eq2.5}) than to directly
deal with ${\mathcal{L}}_n({\mathcal{A}},  {\mathcal{B}})$. In the
model selection procedure introduced in Section~\ref{sec2.3} below, we
may replace ${\mathcal{L}}_n({\mathcal{A}}, {\mathcal{B}})$ by
${\mathcal{L}}_n^\diamond({\mathcal{A}}, {\mathcal{B}})$.

%s2.3 #&#
\subsection{Penalised local likelihood estimation with weighted LASSO
penalties}\label{sec2.3}

In order to conduct the model selection and structure specification, we
define the following penalised local log-likelihood function:
%
%e2.6 #&#
\begin{equation}
\label{eq2.6} {\mathcal Q}_n({\mathcal{A}}, {\mathcal{B}})= {
\mathcal{L}}_n({\mathcal{A}}, {\mathcal{B}})- \sum
_{j=1}^{d_n}p_{nj} \bigl(\| \balpha_j
\| \bigr)- \sum_{j=1}^{d_n}p_{nj}^\ast
\bigl(\| \bbeta_j \| \bigr),
\end{equation}
where $p_{nj}(\cdot)$ and $p_{nj}^\ast(\cdot)$ are two penalty
functions which will be specified later,
\[
\balpha_j = (\alpha_{j1}, \ldots, \alpha_{jn})^{\T}
\quad\mbox{and}\quad \bbeta_j = (\beta_{j1}, \ldots,
\beta_{jn})^{\T},
\]
which correspond to $ [a_j(U_1),  \ldots,  a_j(U_n) ]^{\T}$
and $ [
\dot{a}_j(U_1),  \ldots,  \dot{a}_j(U_n) ]^{\T}$, respectively.
By the quadratic approximation (\ref{eq2.5}), we may approximate
${\mathcal Q}_n({\mathcal{A}},  {\mathcal{B}})$ by ${\mathcal
Q}_n^\diamond({\mathcal{A}},  {\mathcal{B}})$ which is defined through
%
%e2.7 #&#
\begin{equation}
\label{eq2.7} {\mathcal Q}_n^\diamond({\mathcal{A}}, {
\mathcal{B}})= {\mathcal{L}}_n^\diamond({\mathcal{A}}, {
\mathcal{B}})- \sum_{j=1}^{d_n}p_{nj}
\bigl(\| \balpha_j \| \bigr)- \sum_{j=1}^{d_n}p_{nj}^\ast
\bigl(\| \bbeta_j \| \bigr),
\end{equation}
where ${\mathcal{L}}_n^\diamond({\mathcal{A}},  {\mathcal{B}})$ is
defined in (\ref{eq2.5}).

For the penalty functions $p_{nj}(\cdot)$ and $p_{nj}^\ast(\cdot)$ in
(\ref{eq2.6}) and (\ref{eq2.7}), we consider two possible cases: (i)
the adaptive group LASSO penalty, and (ii) the group SCAD penalty. Note
that identifying the constant coefficients in model (\ref{eq2.1}) is
equivalent to identifying the functional coefficients such that either
$\dot{a}_j(U_1) = \cdots= \dot{a}_j(U_n) = 0$ or its deviation
$D_j = 0$, where
\[
D_j = \Biggl\{ \sum_{k=1}^n
\Biggl[ a_j(U_k) - \frac{1}{n} \sum
_{l=1}^n a_j(U_l)
\Biggr]^2 \Biggr\}^{1/2}.
\]
Using the preliminary estimates, we can construct the preliminary
estimator of $D_j$,
\[
\widetilde{D}_j = \Biggl\{ \sum_{k=1}^n
\Biggl[ \widetilde{a}_j(U_k) - \frac{1}{n} \sum
_{l=1}^n \widetilde{a}_j(U_l)
\Biggr]^2 \Biggr\}^{1/2},
\]
where $\widetilde{a}_j(U_k)$ is the $j$th element of $\widetilde{\ba}_k$.

For case (i) of the adaptive group LASSO, we define
%
%e2.8 #&#
\begin{equation}
\label{eq2.8} p_{nj} \bigl(\| \balpha_j \| \bigr)=
\lambda_3\|\widetilde{\balpha}_j \| ^{-\kappa}\|
\balpha_j \|,\qquad p_{nj}^\ast \bigl(\|
\bbeta_j \| \bigr)=\lambda_3^\ast
\widetilde{D}_j^{-\kappa}\| h\bbeta_j \|,
\end{equation}
where $\lambda_3$ and $\lambda_3^\ast$ are two tuning parameters,
$\kappa$ is pre-determined and can be chosen as $1$ or $2$ as in the
literature,
$\widetilde{\balpha}_j
=  [\widetilde{a}_j(U_1),  \ldots,  \widetilde
{a}_j(U_n)
]^{\T}$.

For case (ii) of the group SCAD, we may apply the local linear
approximation to the SCAD penalty function $p_{nj}(\cdot)$ [\citet{ZouLi08}] and then obtain
%
%e2.9 #&#
\begin{equation}
\label{eq2.9} p_{nj} \bigl(\| \balpha_j \| \bigr)\approx
p_{nj} \bigl(\| \widetilde {\balpha}_j \| \bigr) -
\dot{p}_{nj} \bigl(\|\widetilde{\balpha}_j \| \bigr) \|\widetilde{
\balpha}_j \|+\dot{p}_{nj} \bigl(\| \widetilde {
\balpha}_j \| \bigr) \| \balpha_j \|,
\end{equation}
where $p_{nj}(z)\equiv p_{\lambda_4}(z)$ is the SCAD penalty function
with the derivative defined by
%
%e2.10 #&#
\begin{equation}
\label{eq2.10} \dot{p}_{nj}(z)\equiv\dot{p}_{\lambda_4}(z)=
\lambda_4 \biggl[I(z\leq \lambda_4)+
\frac{(a_0\lambda_4-z)_+}{(a_0-1)\lambda}I(z>\lambda _4) \biggr],
\end{equation}
$\lambda_4$ is a tuning parameter and $a_0=3.7$ as suggested in \citet{FanLi01}. Note that the first two terms on the right-hand side of
(\ref
{eq2.9}) do not involve $\|\balpha_j\|$, which motivates us to choose
%
%e2.11 #&#
\begin{equation}
\label{eq2.11} p_{nj} \bigl(\| \balpha_j \| \bigr)=
\dot{p}_{\lambda_4} \bigl(\| \widetilde{\balpha}_j \| \bigr) \|
\balpha_j \|
\end{equation}
with $\dot{p}_{\lambda_4} (\cdot)$ defined in (\ref{eq2.10}). For
$p_{nj}^\ast (\| \bbeta_j \| )$, similar to the corresponding
definition in (\ref{eq2.8}) for case (i), we consider the structure
%
%e2.12 #&#
\begin{equation}
\label{eq2.12} p_{nj}^\ast \bigl(\| \bbeta_j \| \bigr) =
\dot{p}_{\lambda_4^\ast} (\widetilde{D}_j) \| h\bbeta_j
\|,
\end{equation}
where $\dot{p}_{\lambda_4^\ast}(\cdot)$ is defined similar to
$\dot
{p}_{\lambda_4}(\cdot)$ in (\ref{eq2.10}) with $\lambda_4$ replaced by
$\lambda_4^\ast$.

Based on (\ref{eq2.7}) and the above specification of the penalty
functions, we may consider the following two objective functions:
%
%e2.13 #&#
\begin{equation}
\label{eq2.13} {\mathcal Q}_{n}^1({\mathcal{A}}, {
\mathcal{B}}) = {\mathcal{L}}_{n}^\diamond({\mathcal{A}}, {
\mathcal{B}}) - \lambda_3\sum_{j=1}^{d_n}
\|\widetilde{\balpha}_j \|^{-\kappa}\| \balpha _j
\| - \lambda_3^\ast\sum_{j=1}^{d_n}
\widetilde{D}_j^{-\kappa}\| h\bbeta _j \|
\end{equation}
for the adaptive group LASSO penalty; and
%
%e2.14 #&#
\begin{equation}
\label{eq2.14} {\mathcal Q}_{n}^2({\mathcal{A}}, {
\mathcal{B}}) = {\mathcal{L}}_{n}^\diamond({\mathcal{A}}, {
\mathcal{B}}) - \sum_{j=1}^{d_n}
\dot{p}_{\lambda_{4}}\bigl (\| \widetilde{\balpha}_j \| \bigr)\|
\balpha_j \| - \sum_{j=1}^{d_n}
\dot{p}_{\lambda_{4}^\ast} (\widetilde{D}_j)\| h\bbeta _j \|
\end{equation}
for the group SCAD penalty. Note that the penalty terms in (\ref
{eq2.13}) and (\ref{eq2.14}) are the weighted LASSO penalty functions.
In particular, the weights in (\ref{eq2.14}) are determined by the
derivative of the SCAD penalty using the preliminary estimators
$\widetilde{\balpha}_j$ and $\widetilde{D}_j$. The objective functions
in (\ref{eq2.13}) and (\ref{eq2.14}) can be seen, in some sense, as an
extension of that in Bradic, Fan and Wang (\citeyear{BraFanWan11}) from the parametric
linear models to the flexible GSVCMs.

Our model selection and structure specification procedure is based on
maximising the objective function in either (\ref{eq2.13}) or (\ref
{eq2.14}). Let
%
%e2.15 #&#
\begin{eqnarray}
\label{eq2.15} \widehat{\balpha}_j = (\widehat{\alpha}_{j1},
\ldots, \widehat{\alpha}_{jn} )^{\T
} \quad\mbox{and}\quad \widehat{
\bbeta}_j = (\widehat{\beta}_{j1}, \ldots, \widehat{
\beta}_{jn} )^{\T
},
\nonumber
\\[-8pt]
\\[-8pt]
\eqntext{ j=1, \ldots, d_n,}
\end{eqnarray}
be the maximiser of ${\mathcal Q}_{n}^1({\mathcal{A}}, {\mathcal
{B}})$, and
%
%e2.16 #&#
\begin{eqnarray}
\label{eq2.16} \overline{\balpha}_j = (\overline{
\alpha}_{j1}, \ldots, \overline{\alpha}_{jn} )^{\T}
\quad\mbox{and} \quad\overline{\bbeta}_j = (\overline{\beta}_{j1},
\ldots, \overline{\beta}_{jn} )^{\T
},
\nonumber
\\[-8pt]
\\[-8pt]
\eqntext{j=1, \ldots,
d_n,}
\end{eqnarray}
be the maximiser of ${\mathcal Q}_{n}^2({\mathcal{A}},  {\mathcal
{B}})$. The asymptotic theorems and remarks in Section~\ref{sec3} show
that the estimators defined in (\ref{eq2.15}) and (\ref{eq2.16}) equal
to the biased oracle estimators of the functional coefficients (see
Section~\ref{sec3} for the definition) with probability approaching one.

%s2.4 #&#
\subsection{Estimation of the constant coefficients}\label{sec2.4}

We next discuss how to estimate the constant coefficients in GSVCMs. By
choosing the penalty function as the adaptive group LASSO (or the group
SCAD) penalty, we would expect $\|\widehat{\balpha}_j \| = 0$ (or $\|
\overline{\balpha}_j \| = 0$) when $a_j(\cdot) = 0$, and $\|\widehat
{\bbeta}_j \| = 0$ (or $\|\overline{\bbeta}_j \| = 0$) when
$a_j(\cdot
)$ is a constant. Hence, our model selection and structure
specification procedure works as follows: if $\|\widehat{\balpha}_j \|
= 0$ (or $\|\overline{\balpha}_j \| = 0$), the corresponding variable
$x_j$ is not significant and should be removed from
the model; if $\|\widehat{\bbeta}_j \| = 0$ (or $\|\overline{\bbeta}_j
\| = 0$), the functional coefficient of $a_j(\cdot)$ is a constant
which is denoted by $c_j$ and can be estimated by
%
%e2.17 #&#
\begin{equation}\quad
\label{eq2.17} \widehat{c}_j = n^{-1} \sum
_{i=1}^n \widehat{\alpha}_{ji}\quad \mbox{or}\quad
\overline{c}_j = n^{-1} \sum_{i=1}^n
\overline{\alpha }_{ji},\qquad j=s_{n1}+1,\ldots,s_{n2}.
\end{equation}
Then the semi-varying coefficient modelling structure is finally specified.

%%%%%%%%%%%%%%%%%%%%%%

%s3 #&#
\section{Asymptotic theory}
\label{sec3}

In this section, we present the asymptotic properties of the model
selection and structure specification procedure introduced in Section~\ref{sec2}. Recall that
\[
\ba_{k0}= \bigl[a_1(U_k), \ldots,
a_{d_n}(U_k) \bigr]^{\T} \quad\mbox{and}\quad
\bb_{k0}= \bigl[\dot{a}_1(U_k), \ldots, \dot
{a}_{d_n}(U_k) \bigr]^{\T}
\]
for $k=1,\ldots,n$. We start with the uniform consistency results for
their penalised local maximum
likelihood estimators
\[
\widetilde{\ba}_k= \bigl[\widetilde{a}_1(U_k),
\ldots, \widetilde {a}_{d_n}(U_k) \bigr]^{\T}\quad
\mbox{and}\quad \widetilde{\bb}_k= \bigl[\widetilde{\dot{a}}_1(U_k),
\ldots, \widetilde{\dot {a}}_{d_n}(U_k)
\bigr]^{\T},
\]
which are the maximisers of the objective function in (\ref{eq2.3}). In
the sequel, we let $\alpha_n \propto\beta_n$ denote $c_1\beta_n\leq
\alpha_n\leq c_2\beta_n$ when $n$ is sufficiently large, where
$0<c_1\leq c_2<\infty$.

\begin{prop}\label{pr3.1}
Suppose that Assumptions \ref{ass1}--\ref{ass4} in
Appendix \ref{app.A} are satisfied.

\textup{(i)} If the moment condition (\ref{eqA.1}) and Assumption~\ref{ass5}
are satisfied with $d_n\propto n^{\tau_1}$, $0\leq\tau_1<\infty$,
we have
%
%e3.1 #&#
\begin{equation}
\label{eq3.1} \max_{1\leq k\leq n}\Vert\widetilde{\ba}_k-
\ba_{k0}\Vert+ \max_{1\leq
k\leq n}\bigl\Vert h(\widetilde{
\bb}_k-\bb_{k0})\bigr\Vert=O_P(\sqrt
{s_{n2}}\lambda_1),
\end{equation}
where $s_{n2}$ is the number of the nonzero functional coefficients

\textup{(ii)} If the moment condition (\ref{eqA.2}) and Assumption
\ref{ass5prime} are satisfied with $d_n\propto\exp \{(nh)^{\tau
_2} \}$, then (\ref{eq3.1}) also holds, where $0\leq\tau
_2<1-\tau_3$
with $0<\tau_3<1$.
\end{prop}

\begin{rem}\label{re3.1} The above proposition indicates that the
preliminary estimators $\widetilde{\ba}_k$ and $\widetilde{\bb}_k$ are
uniformly consistent, as Assumption~\ref{ass3} in Appendix \ref{app.A}
guarantees that the maximal distance between two consecutive index
variables $U_i$ is only with the order $O_P(\log n/n)$ [cf., \citet{Jan87}] and the observed values of $U$ can be sufficiently dense on the
compact support $[0,1]$. The uniform convergence rate in (\ref{eq3.1})
depends on $s_{n2}$, the number of the nonzero functional
coefficients, and the tuning parameter $\lambda_1$. In Assumptions \ref{ass5}
and \ref{ass5prime}, we impose some condition on the relationship between
$\lambda_1$ and the well-known uniform convergence rate $ (\frac
{\log h^{-1}}{nh} )^{1/2}$, and assume that $\lambda_1\propto
\lambda
_2$. As a consequence, the influence of $ (\frac{\log
h^{-1}}{nh}
)^{1/2}$ and $\lambda_2$ would be dominated by that of $\lambda_1$.
Although the dimension of potential covariates in our model can be
larger than the sample size and diverge at an exponential rate,
$s_{n2}$ is not allowed to diverge too fast in order to guarantee the
consistency of the preliminary estimators of the functional
coefficients. The condition $s_{n2}\lambda_1^2h^{-2}=o(1)$ in
Assumptions \ref{ass5} and \ref{ass5prime} indicates that $s_{n2}$ is allowed to be
divergent at a slow polynomial rate of $n$. It is also interesting to
find from the comparison between (\ref{eqA.1}) and (\ref{eqA.2}) in
Appendix \ref{app.A} that the required moment condition when $d_n$ diverges at a
polynomial rate is weaker than that when $d_n$ diverges at an
exponential rate.
\end{rem}

\begin{rem}\label{re3.2} Note that in the penalised local
log-likelihood estimation method in Section~\ref{sec2.1}, we do not use
the group LASSO or SCAD penalty function. Although Proposition~\ref{pr3.1}
establishes the uniform consistency for the preliminary estimators of
the functional coefficients and their derivatives, the shrinkage
estimation method in Section~\ref{sec2.1} does not have the well-known
sure screening property [Fan, Ma and Dai (\citeyear{FanMaDai14}); Liu, Li and Wu
(\citeyear{LiuLiWu14})]. However, under some further conditions, the uniform convergence
rate in (\ref{eq3.1}) would be sufficient for us to prove Theorems \ref{th3.1}
and \ref{th3.2} below.

Let
\[
\overline{\mathcal{A}}_n= \bigl(\overline{\ba}_1^{\T},
\ldots, \overline{\ba} _n^{\T} \bigr)^{\T}\quad
\mbox{and}\quad \overline{\mathcal{B}}_n= \bigl(\overline{
\bb}_1^{\T}, \ldots, \overline{\bb} _n^{\T}
\bigr)^{\T},
\]
where $\overline{\ba}_k = (\overline{\alpha}_{1k}, \ldots,
\overline
{\alpha}_{d_nk})^{\T}$ and $\overline{\bb}_k = (\overline{\beta}_{1k},
\ldots, \overline{\beta}_{d_nk})^{\T}$. Let $\ba_k^o$ be any
$d_n$-dimensional vector with the last $(d_n-s_{n2})$ elements being
zeros, and $\bb_k^o$ be any $d_n$-dimensional vector with the last
$(d_n-s_{n1})$ elements being zeros. Denote
\[
{\mathcal{A}}^{o}= \bigl[\bigl(\ba_1^{o}
\bigr)^{\T}, \ldots, \bigl(\ba_n^{o}
\bigr)^{\T
} \bigr]^{\T} \quad\mbox{and}\quad {\mathcal{B}}^{o}=
\bigl[\bigl(\bb_1^{o}\bigr)^{\T}, \ldots, \bigl(
\bb_n^{o}\bigr)^{\T} \bigr]^{\T},
\]
and then define the biased oracle estimators
\[
\overline{\mathcal{A}}_n^{bo}= \bigl[\bigl(\overline{
\ba}_1^{bo}\bigr)^{\T}, \ldots, \bigl(\overline{
\ba}_n^{bo}\bigr)^{\T} \bigr]^{\T}\quad
\mbox{and}\quad \overline {\mathcal {B}}_n^{bo}= \bigl[\bigl(
\overline{\bb}_1^{bo}\bigr)^{\T}, \ldots, \bigl(
\overline {\bb} _n^{bo}\bigr)^{\T}
\bigr]^{\T},
\]
which maximise the objective function ${\mathcal Q}_{n}^2({\mathcal
{A}}^o, {\mathcal{B}}^o)$ when the penalty function is the SCAD
penalty. The following theorem gives the relation between the penalised
estimators which maximise the objective function (\ref{eq2.14}) and the
corresponding biased oracle estimators when the SCAD penalty is used.
The result for the case of the adaptive group LASSO penalty is similar
and will be discussed in Remark~\ref{re3.3} below.
\end{rem}

\begin{theorem}\label{th3.1}
Suppose that the conditions in
Proposition~\ref{pr3.1}\textup{(i)} are satisfied. When the penalty functions are
defined in (\ref{eq2.11}) and (\ref{eq2.12}), and Assumption~\ref{ass6} in
Appendix \ref{app.A} is satisfied, with probability approaching one,
the maximiser of the objective function ${\mathcal Q}_n^2(\cdot,\cdot)$
defined in (\ref{eq2.14}), $(\overline{\mathcal{A}}_n, \overline
{\mathcal{B}}_n)$, exists and equals to
$(\overline{\mathcal{A}}_n^{bo}, \overline{\mathcal{B}}_n^{bo})$.
Furthermore,
%
%e3.2 #&#
\begin{equation}
\label{eq3.2} \frac{1}{n} \bigl\Vert\overline{\mathcal A}_n^{bo}-{
\mathcal A}_0\bigr \Vert^2=\frac{s_{n2}}{nh},\qquad
\frac{1}{n} \bigl\Vert\overline {\mathcal B}_n^{bo}-{
\mathcal B}_0 \bigr\Vert^2=\frac{s_{n2}}{nh^3},
\end{equation}
where
\[
{\mathcal A}_0 = \bigl(\ba_{10}^{\T}, \ldots,
\ba_{n0}^{\T}\bigr)^{\T},\qquad {\mathcal B}_0
= \bigl(\bb_{10}^{\T}, \ldots, \bb_{n0}^{\T}
\bigr)^{\T}.
\]
\end{theorem}

\begin{rem}\label{re3.3}
Given the moment condition (\ref{eqA.2})
and Assumption~\ref{ass5prime} in Appendix \ref{app.A} with $d_n\propto
\exp
 \{(nh)^{\tau_2} \}$, the above result still holds. It can be
proved by using Proposition~\ref{pr3.1}(ii) and strengthening (\ref{eqA.3}) in
Assumption~\ref{ass6} to
\[
h^{-1/2} \biggl[ \biggl(\frac{\log h^{-1}}{nh} \biggr)^{\tau_3/2}\sqrt
{nh}+s_{n2}^{1/2}(1+\lambda_1\sqrt{nh}) \biggr]=o(
\lambda_4),
\]
where $\tau_3$ is defined in Proposition~\ref{pr3.1}(ii). Noting that the
left-hand side is controlled by $\lambda_1\sqrt{ns_{n2}}$, the above
condition can be simplified to $\lambda_1\sqrt{ns_{n2}}=o(\lambda_4)$.
Theorem~\ref{th3.1} suggests, using the proposed model selection procedure, the
zero coefficients can be estimated exactly as zeros, and the
derivatives of the
constant coefficients can also be estimated exactly as zeros, which indicates
that the \textit{sparsity} property holds for the proposed model selection
procedure. Hence, our theorem complements some existing ultra-high
dimensional sparsity results such as those derived by Bradic, Fan and Wang (\citeyear{BraFanWan11}), \citet{FanLv11} and \citet{Lia12}. Furthermore, for the
penalty functions defined in (\ref{eq2.11}) and (\ref{eq2.12}), by
Proposition~\ref{pr3.1} and Assumption~\ref{ass6}, we may show that properties (i)--(iv)
for the folded concave penalty function introduced by Fan, Ma and Dai
(\citeyear{FanMaDai14}) are satisfied with probability approaching one. Hence, Theorem~\ref{th3.1} can also be seen, in some sense, as a generalisation of Theorem~1
in Fan, Ma and Dai (\citeyear{FanMaDai14}). When the adaptive group LASSO penalty is
used, by modifying the conditions in Assumption~\ref{ass6}(i), we may show that
the above sparsity result still holds and (\ref{eq3.2}) is satisfied by
replacing $\overline{\mathcal{A}}_n^{bo}$ and $\overline{\mathcal
{B}}_n^{bo}$ by $\widehat{\mathcal{A}}_n^{bo}$ and $\widehat
{\mathcal
{B}}_n^{bo}$, the biased oracle estimators with the adaptive group
LASSO penalty.
\end{rem}

We next study the oracle property for the penalised estimators of the nonzero
functional coefficients and constant coefficients. Let
$a_{j}^{uo}(U_k)$, $j=1, \ldots, s_{n1}$,
$k=1,\ldots,n$, be the (unbiased) oracle estimator of $a_j(U_k)$, and
$c_{j}^{uo}$,
$j=s_{n1}+1, \ldots, s_{n2}$, be the (unbiased)\vspace*{1pt} oracle estimator of
the constant
coefficient $c_j$. The (unbiased) oracle estimators are obtained by the standard
estimation procedure for the GSVCMs, that is, the maximisation of the objective
function ${\mathcal{L}}_{n}^\diamond({\mathcal{A}}^o, {\mathcal
{B}}^o)$ with
respect to ${\mathcal{A}}^o$ and ${\mathcal{B}}^o$ [the penalty terms in
(\ref{eq2.13}) and (\ref{eq2.14}) are ignored] and the application of
(\ref{eq2.17}) under the assumption that we know $a_j(\cdot) \equiv
0$ when
$j = s_{n2} + 1, \ldots, d_n$ and $a_j(\cdot)\equiv c_j$ when
$j = s_{n1}+1,\ldots,s_{n2} $. In the following theorem, we only
consider the case that the penalty functions are defined in (\ref
{eq2.11}) and (\ref{eq2.12}) to save the space. Let
\[
\overline{\mathbf{D}}_n = \Bigl( \max_{1 \leq k \leq n} \bigl|
\overline{a}_{1}(U_k) - a_{1}^{uo}(U_k)
\bigr|, \ldots, \max_{1 \leq k \leq n} \bigl| \overline{a}_{s_{n1}}(U_k)
- a_{s_{n1}}^{uo}(U_k) \bigr| \Bigr)^{\T},
\]
where $\overline{a}_{j}(U_k)=\overline{\alpha}_{jk}$ is defined in
(\ref
{eq2.16}), and
\[
{\bC}_n^{uo} = \bigl( c_{s_{n1}+1}^{uo},
\ldots, c_{s_{n2}}^{uo} \bigr)^{\T},\qquad \overline{
\bC}_n = ( \overline{c}_{s_{n1}+1}, \ldots, \overline{c}_{s_{n2}}
)^{\T},
\]
where $\overline{c}_j$ is defined in (\ref{eq2.17}).

\begin{theorem}\label{th3.2}
Suppose that the
conditions of
Theorem~\ref{th3.1} are satisfied. For any $s_{n1}$-dimensional vector
${\mathbf
{B}}_n$ with
$\|{\mathbf{B}}_n \| = 1$, we have
%
%e3.3 #&#
\begin{equation}
\label{eq3.3} \sqrt{nh} {\mathbf{B}}_n^{\T} \overline{
\mathbf{D}}_n = o_P(1);
\end{equation}
and for any $(s_{n2}-s_{n1})$-dimensional vector ${\mathbf{A}}_n$ with
$\|{\mathbf{A}}_n\| = 1$, we have
%
%e3.4 #&#
\begin{equation}
\label{eq3.4} \sqrt{n} {\mathbf{A}}_n^{\T} \bigl(
\overline{\bC}_n - {\bC}_n^{uo} \bigr) =
o_P(1).
\end{equation}
\end{theorem}

\begin{rem}\label{re3.4} Theorem~\ref{th3.2} indicates that the penalised
likelihood estimators of the nonzero functional coefficients and
constant coefficients have the same asymptotic distribution as the
corresponding oracle estimators, and thus the oracle property holds. As
discussed in Remark~\ref{re3.3}, by strengthening the moment conditions, we can
also show that the above oracle property holds when $d_n\propto\exp
\{(nh)^{\tau_2} \}$. Following the arguments in \citet{ZhaPen10} and Li, Ke and Zhang (\citeyear{LiKeZha13}), we can easily establish the
asymptotic normality of $\overline{a}_j(\cdot)$, $j=1, \ldots,
s_{n1}$, and $\overline{c}_j$, $j=s_{n1}+1, \ldots, s_{n2}$.
\end{rem}
%%%%%%%%%%%%%%%%%%%%%%

%s4 #&#
\section{Computational algorithm and selection of tuning parameters}
\label{sec4}

In this section, we introduce a computational algorithm to maximise
${\mathcal Q}_n^1({\mathcal{A}}, {\mathcal{B}})$ and ${\mathcal
Q}_n^2({\mathcal{A}}, {\mathcal{B}})$ defined in Section~\ref{sec2.3}
and discuss how to
choose the tuning parameters involved in the proposed penalised likelihood
method.

%s4.1 #&#
\subsection{Computational algorithm}\label{sec4.1}

We first re-arrange the quadratic objective function ${\mathcal
{L}}_n^{\diamond}({\mathcal{A}}, {\mathcal{B}})$ in order to make it
have the standard form
when using the penalised estimation method. Let
\[
\btheta= \bigl(\balpha_1^{\T}, \ldots,
\balpha_{d_n}^{\T}, h\bbeta_1^{\T},
\ldots,h\bbeta_{d_n}^{\T} \bigr)^{\T}
\]
and define the transformation matrix
\[
\bT= (I_n \otimes e_{1,2d_n}, \ldots, I_n \otimes
e_{d_n,2d_n}, I_n \otimes e_{d_n+1,2d_n}, \ldots,
I_n \otimes e_{2d_n,2d_n} )^{\T},
\]
where $e_{k,d}$ is a $d$-dimensional unit vector with the $k$th
component being 1 and $I_n$ is an $n\times n$ identity matrix. With the
above notation, it is easy to show that $\btheta= {\bT} {\mathcal
V}_n({\mathcal{A}}, h{\mathcal{B}})$,
where ${\mathcal V}_n({\mathcal{A}}, h{\mathcal{B}})$ is defined as in
Section~\ref{sec2.2}. Let $\tilde{\btheta}$ be defined as $\btheta$ but
with $\calA$ and $\calB$ replaced by $\tilde{\calA}$ and $\tilde
{\calB
}$, respectively, and
\[
\bH^2=\bH^{\T}\bH=- \bT\ddot{\mathcal{L}}_n(
\widetilde{\mathcal{A}}_n, \widetilde{\mathcal{B}}_n)
\bT^{\T}, \qquad\tilde{\eeta} = \bH\tilde{\btheta} + \bigl({\bH}^{-1}
\bigr)^{\T} \bT\dot {\mathcal {L}}_n(\widetilde{
\mathcal{A}}_n, \widetilde{\mathcal{B}}_n).
\]
We define a quadratic objective function
\[
\calL_{n}^{\ast} ({\mathcal{A}}, {\mathcal{B}}) =-
\tfrac{1}{2} (\tilde{\eeta}- \bH\btheta)^{\T} (\tilde{\eeta}- \bH
\btheta).
\]
Given the initial estimator
${\mathcal V}_n(\widetilde{\mathcal{A}}_n, h\widetilde{\mathcal
{B}}_n)$, it is
easy to see the difference between
${\mathcal{L}}_n^{\diamond}({\mathcal{A}}, {\mathcal{B}})$ and
$\calL_{n}^{\ast} ({\mathcal{A}}, {\mathcal{B}})$ is a constant.\vspace*{1pt}
Therefore,
the maximiser of ${\mathcal Q}_n^1({\mathcal{A}}, {\mathcal{B}})$ or
${\mathcal Q}_n^2({\mathcal{A}}, {\mathcal{B}})$ is the minimiser
of the following target function:
%
%e4.1 #&#
\begin{equation}
\label{eq4.2} \calO(\btheta) \equiv \frac{1}{2} (\tilde{\eeta}- \bH
\btheta)^{\T} (\tilde{\eeta}- \bH\btheta) + \sum
_{j=1}^{d_n} \tau_{1j} \|
\balpha_j \|+ \sum_{j=1}^{d_n}
\tau_{2j} \| h\bbeta_j \|,
\end{equation}
where
\[
\tau_{1j} = \lambda_3 \|\widetilde{\balpha}_j
\|^{-\kappa}\quad \mbox{and}\quad \tau_{2j} = \lambda_3^\ast
\widetilde{D}_j^{-\kappa}
\]
for ${\mathcal Q}_n^1({\mathcal{A}}, {\mathcal{B}})$; and
\[
\tau_{1j} = \dot{p}_{\lambda_{4}} \bigl(\| \widetilde{
\balpha}_j \| \bigr) \quad\mbox{and}\quad \tau_{2j} =
\dot{p}_{\lambda_{4}^\ast} (\widetilde{D}_j)
\]
for ${\mathcal Q}_n^2({\mathcal{A}},  {\mathcal{B}})$. As a direct
consequence of the Karush--Kuhn--Tucker conditions, we have that a
necessary and sufficient condition for $\btheta$ to be a minimiser of
$\calO(\btheta)$ is, for $j=1,\ldots,d_n$,
\[
\cases{ %
-\mathrm{H}_j^{\T}(
\tilde{\eeta}- \bH\btheta) + \tau_{1j} \| \balpha _j \|
^{-1} \balpha_j = {\mathbf0}_n, &\quad $\forall
\balpha_j \neq{\mathbf0}_n,$
\vspace*{2pt}\cr
\bigl\|\mathrm{H}_j^{\T}(\tilde{\eeta}- \bH\btheta)\bigr\| <
\tau_{1j}, & \quad $\forall \balpha_j = {\mathbf0}_n,$
\vspace*{2pt}\cr
-\mathrm{H}_{j+d_n}^{\T}(\tilde{\eeta}- \bH\btheta) +
\tau_{2j} \| \bbeta _j \|^{-1}
\bbeta_j= {\mathbf0}_n, &\quad $\forall \bbeta_j
\neq{\mathbf0}_n,$
\vspace*{2pt}\cr
\bigl\|\mathrm{H}_{j+d_n}^{\T}(\tilde{\eeta}- \bH\btheta)\bigr\| <
\tau_{2j}, & \quad$\forall \bbeta_j = {\mathbf0}_n,$}
\]
where $\mathrm{H}_j$ is the matrix consisting of the $((j-1)n+1)$th to
the $(jn)$th column of $\bH$ and $\0_n$ is an $n$-dimensional vector
with each component being $0$. Hence, for $j=1,  \ldots,  d_n$, we
have $\balpha_j={\mathbf0}_n$ if $\|\mathrm{H}_j^{\T}(\tilde{\eeta}- \bH
\btheta_{-j})\| < \tau_{1j}$, otherwise
\[
\balpha_j= \bigl(\mathrm{H}_j^{\T}
\mathrm{H}_j + \tau_{1j} \|\balpha_j
\|^{-1} I_n \bigr)^{-1} \mathrm{H}_j^{\T}(
\tilde{\eeta}- \bH\btheta_{-j});
\]
and $\bbeta_j ={\mathbf0}_n$ if $\|\mathrm{H}_{j+d_n}^{\T}(\tilde{\eeta}-
\bH
\btheta_{-(j+d_n)})\| < \tau_{2j}$, otherwise
\[
\bbeta_j = \bigl(h\mathrm{H}_{j+d_n}^{\T}
\mathrm{H}_{j+d_n} + \tau_{2j} \| \bbeta_j
\|^{-1} I_n \bigr)^{-1} \mathrm{H}_{j+d_n}^{\T}(
\tilde{\eeta}- \bH\btheta_{-(j+d_n)}),
\]
where
\begin{eqnarray*}
\btheta_{-j} &=& \bigl(\balpha_1^{\T}, \ldots,
\balpha_{j-1}^{\T}, \0_n^{\T},
\balpha_{j+1}^{\T}, \ldots, \balpha_{d_n}^{\T},
h\bbeta_1^{\T}, \ldots, h\bbeta_{d_n}^{\T}
\bigr)^{\T},
\\
\btheta_{-(j+d_n)} &=& \bigl(\balpha_1^{\T}, \ldots,
\balpha_{d_n}^{\T}, h\bbeta_1^{\T}, \ldots,
h\bbeta_{j-1}^{\T}, \0_n^{\T}, h\bbeta
_{j+1}^{\T}, \ldots, h\bbeta_{d_n}^{\T}
\bigr)^{\T}.
\end{eqnarray*}

This leads to the following iterative algorithm to obtain the
minimisers of $\calO(\btheta)$.
\begin{longlist}[\textit{Step} 1.]
\item[\textit{Step} 1.] Start with $\balpha_j^{(0)} = \widetilde
{\balpha
}_j$ and
$\bbeta_j^{(0)} = \widetilde{\bbeta}_j$, $j=1,  \ldots,\break  d_n$, where
$\widetilde{\balpha}_j$ and $\widetilde{\bbeta}_j$ are the preliminary
estimates of the functional coefficients $ [a_j(U_1), \ldots,
a_j(U_n) ]^{\T}$ and their derivatives $ [\dot
{a}_j(U_1),\ldots
, \dot{a}_j(U_n) ]^{\T}$, respectively, which are introduced in
Section~\ref{sec2.1}.

\item[\textit{Step} 2.] For $j=1,  \ldots,  d_n$, let $\balpha_j^{(k)}$
and $\bbeta_j^{(k)}$ be the results after the $k$th iteration. Update
$\balpha_j^{(k)}$ and $\bbeta_j^{(k)}$ in the $(k+1)$th iteration as
follows: for $j=1,  \ldots,  d_n$, $\balpha_j^{(k+1)} = {\mathbf0}_n$ if
$\|\mathrm{H}_j^{\T}(\tilde{\eeta}- \bH\btheta_{-j}^{(k)})\| < \tau
_{1j}^{(k)}$, otherwise
\[
\balpha_j^{(k+1)}= \bigl(\mathrm{H}_j^{\T}
\mathrm{H}_j + \tau_{1j}^{(k)} \bigl\|
\balpha_j^{(k)}\bigr\|^{-1} I_n
\bigr)^{-1} \mathrm{H}_j^{\T}\bigl(\tilde{\eeta}-
\bH\btheta_{-j}^{(k)}\bigr);
\]
and $\bbeta_j^{(k+1)} = {\mathbf0}_n$ if $\|\mathrm{H}_{j+d_n}^{\T}(\tilde
{\eeta}- \bH\btheta_{-(j+d_n)}^{(k)})\| < \tau_{2j}^{(k)}$, otherwise
\[
\bbeta_j^{(k+1)} = \bigl( h\mathrm{H}_{j+d_n}^{\T}
\mathrm{H}_{j+d_n} + \tau_{2j}^{(k)} \bigl\|\bbeta
_j^{(k)}\bigr\|^{-1} I_n
\bigr)^{-1} \mathrm{H}_{j+d_n}^{\T}\bigl(\tilde{\eeta}-
\bH\btheta_{-(j+d_n)}^{(k)}\bigr);
\]
where $\tau_{1j}^{(k)}$ is defined as $\tau_{1j}$ in (\ref{eq4.2}) but
with $\widetilde{\balpha}_j$ replaced by ${\balpha}_j^{(k)}$, $\tau
_{2j}^{(k)}$ is defined as $\tau_{2j}$ in (\ref{eq4.2}) but with
$\widetilde{D}_j$ replaced by ${D}_j^{(k)}$,
\begin{eqnarray*}
D_j^{(k)}&=& \Biggl\{\sum_{s=1}^n
\Biggl[a_j^{(k)}(U_s) - \frac
{1}{n}\sum
_{l=1}^na_j^{(k)}(U_l)
\Biggr]^2 \Biggr\}^{1/2},
\\
\btheta_{-j}^{(k)}&= &\bigl[\bigl(\balpha_1^{(k+1)}
\bigr)^{\T}, \ldots, \bigl(\balpha _{j-1}^{(k+1)}
\bigr)^{\T}, \0_n^{\T}, \bigl(\balpha_{j+1}^{(k)}
\bigr)^{\T}, \ldots, \bigl(\balpha_{d_n}^{(k)}
\bigr)^{\T},
\\
& &{}\bigl(h\bbeta_1^{(k)}\bigr)^{\T}, \ldots,
\bigl(h\bbeta_{d_n}^{(k)}\bigr)^{\T
}
\bigr]^{\T}\quad\mbox{and}
\\
\btheta_{-(j+d_n)}^{(k)}&=&\bigl[\bigl(\balpha_1^{(k+1)}
\bigr)^{\T}, \ldots, \bigl(\balpha_{d_n}^{(k+1)}
\bigr)^{\T}, \bigl(h\bbeta_1^{(k+1)}
\bigr)^{\T}, \ldots, \bigl(h\bbeta _{j-1}^{(k+1)}
\bigr)^{\T}, \0_n^{\T},
\\
&& \bigl(h\bbeta_{j+1}^{(k)}\bigr)^{\T}, \ldots,
\bigl(h\bbeta_{d_n}^{(k)}\bigr)^{\T
}
\bigr]^{\T}.
\end{eqnarray*}
Furthermore, if $\|\balpha_j^{(k)}\| = {\mathbf0}_n$ and $\|\mathrm{H}_j^{\T
}(\tilde{\eeta}- \bH\btheta_{-j}^{(k)})\| > \tau_{1j}^{(k)}$, we set
\[
\balpha_j^{(k+1)}= \bigl(\mathrm{H}_j^{\T}
\mathrm{H}_j + \bigl(\tau _{1j}^{(k)}/
\Delta^{(k)}_{\balpha}\bigr) I_n \bigr)^{-1}
\mathrm{H}_j^{\T
}\bigl(\tilde{\eeta}- \bH
\btheta_{-j}^{(k)}\bigr)
\]
with $\Delta^{(k)}_{\balpha}= \min \{\|\balpha_l^{(k)}\|: \|
\balpha
_l^{(k)}\| \neq0,  l = 1, \ldots, d_n \}$. If $\|\bbeta
_j^{(k)}\| ={\mathbf0}_n$ and\break   $\|\mathrm{H}_{j+d_n}^{\T}(\eeta- \bH\btheta
_{-(j+d_n)}^{(k)})\| > \tau_{2j}^{(k)}$, we set
\[
\bbeta_j^{(k+1)}= \bigl(h\mathrm{H}_{j+d_n}^{\T}
\mathrm{H}_{j+d_n} + \bigl(\tau _{2j}^{(k)}/
\Delta_{\bbeta}^{(k)}\bigr) I_n \bigr)^{-1}
\mathrm{H}_{j+d_n}^{\T
}\bigl(\tilde{\eeta}- \bH
\btheta_{-(j+d_n)}^{(k)}\bigr)
\]
with $\Delta_{\bbeta}^{(k)}=\min \{\|\bbeta_l^{(k)}\|: \|
\bbeta
_l^{(k)}\| \neq0,  l = 1,  \ldots,  d_n
 \}$.

\item[\textit{Step} 3.] If $\sum_{j=1}^{d_n} (\|\balpha
_j^{(k)} - \balpha_j^{(k+1)}\| + h \|\bbeta_j^{(k)} - \bbeta
_j^{(k+1)}\|
 )$ is smaller than a chosen threshold, we stop the iteration, and
$ (\balpha_j^{(k+1)},  \bbeta_j^{(k+1)} )$, $j=1,  \ldots,
 d_n$, are the minimisers of $\calO(\btheta)$.
\end{longlist}

The simulation studies in Section~\ref{sec5} below will show that the
above iterative procedure works reasonably well in the finite sample
cases. The simulation studies are conducted by a small computer cluster
which contains 64 CPUs while the real data analysis results are
obtained by a single PC within one day.

%%%%%%%%%%%%%%%%%%%%%%

%s4.2 #&#
\subsection{Selection of the tuning parameters}\label{sec4.2}

The tuning parameters involved in the proposed model selection and
structure specification procedure play a very important role. We next
discuss how to choose these tuning parameters. First, for the
preliminary estimates, the tuning parameters $\lambda_1$ and $\lambda
_2$ are selected through BIC, and the bandwidth is set to be $h =
0.75[(\log d_n) / n]^{0.2}$. The reason for not using a data-driven
method to select the bandwidth $h$ is to reduce the computational cost.
Also the preliminary estimation is not very sensitive to the choice of
the bandwidth. Then, for the model selection and specification
procedure based on ${\mathcal Q}_n^1({\mathcal{A}},  {\mathcal{B}})$
or ${\mathcal Q}_n^2({\mathcal{A}},  {\mathcal{B}})$, the tuning
parameters $\lambda_3$ and $\lambda_3^{\ast}$ or $\lambda_4$ and
$\lambda_4^{\ast}$ are selected by the generalised information
criterion (GIC) proposed by \citet{FanTan13}. We next briefly
introduce the GIC method.

As the models concerned involve both unknown constant parameters and
unknown functional parameters, to use GIC, we first need to figure out
how many unknown constant parameters an unknown functional parameter
amounts to. Cheng, Zhang and Chen (\citeyear{CheZhaChe09}) suggest that an unknown
functional parameter would amount to $1.028571 h^{-1}$ unknown constant
parameters when Epanechnikov kernel is used. Hence, we construct the
GIC for model (\ref{eq2.1}) as
\begin{eqnarray*}
\GIC\bigl(\lambda, \lambda^{\ast}\bigr) &=& -2\sum
_{i=1}^{n}\ell\bigl(\hat{m}(U_i,
X_i), y_{i}\bigr)
\\
&&+{}2\mathrm{ln\bigl\{ln}(n)\bigr\}\mathrm{ln}\bigl(1.028571 d_{n}h^{-1}
\bigr) \bigl(k_{1}+1.028571 k_{2}h^{-1}\bigr),
\end{eqnarray*}
where $\hat{m}(U_i, X_i)$ is defined as $m(U_i, X_i)$ with all unknowns
being replaced by their estimators obtained based on the tuning
parameters $\lambda_3$ and $\lambda_3^{\ast}$ (or $\lambda_4$ and
$\lambda_4^{\ast}$), $k_{1}$ is the number of significant covariates
with constant
coefficients obtained based on the given pair of tuning parameters, and
$k_{2}$ is the number of significant covariates with
functional coefficients obtained based on the given pair of tuning
parameters. For the maximisation of ${\mathcal Q}_n^1({\mathcal{A}},  {\mathcal{B}})$, the minimiser of $\GIC(\lambda_3, \lambda_3^{\ast})$
is the selected $\lambda_3$ and $\lambda_3^{\ast}$, while for the
maximisation of ${\mathcal Q}_n^2({\mathcal{A}},  {\mathcal{B}})$, the
minimiser of $\GIC(\lambda_4, \lambda_4^{\ast})$ is the selected
$\lambda_4$ and~$\lambda_4^{\ast}$.

%%%%%%%%%%%%%%%%%%%

%s5 #&#
\section{Simulation studies}\label{sec5}

In this section, we give three simulated examples to examine the
accuracy of the proposed model selection, structure specification and
estimation procedure, as well as the oracle property of the proposed
estimators. Throughout this section, we call the procedure based on
(\ref{eq2.13}) the adaptive group LASSO method and the procedure based
on (\ref{eq2.14}) the group SCAD method. For the adaptive group LASSO
method, the pre-determined parameter $\kappa$ is chosen to be~1. For
the group SCAD method, the SCAD penalty is defined through its
derivative as in (\ref{eq2.10}). The kernel function used in this
section is taken to be the Epanechnikov kernel $K(t) = 0.75 (1-t^2)_+$.
The bandwidth and other tuning parameters are selected by the approach
described in Section~\ref{sec4.2}.

We will start with a simulated example on a semi-varying coefficient
Poisson regression model, then an example on varying coefficient models
and finally an example on a varying coefficient logistic regression
model. In Example~\ref{ex5.1}, we will compare the performance of the proposed
adaptive group LASSO and group SCAD methods on model selection,
structure specification and estimation, and find that the group SCAD
method gives slightly better finite sample performance under all
simulation settings.
Thus, we will call the group SCAD method ``our method'' in the following
two examples and only compare it with some existing methods. In Example~\ref{ex5.2}, we will compare our method with the KLASSO proposed in \citet{WanXia09} based on varying coefficient models. In Example~\ref{ex5.3}, we will
compare our method with the method proposed in \citet{Lia12}. The
simulation results of the KLASSO and Lian's method in Tables~\ref{tab3}--\ref{tab5} are the original results reported in \citet{WanXia09} and \citet{Lia12}, respectively. From the simulation results, we
will find that our method outperforms the existing ones.

\begin{exa}\label{ex5.1}
We generate a sample from a Poisson
regression model as follows: first independently generate $x_{ij}$,
$i=1,\ldots, n$, $j=1, \ldots,d_n$, from the standard normal
distribution $\mathsf{ N}(0,  1)$, and $U_i$, $i=1, \ldots, n$, from
uniform distribution $\mathsf{ U}[0,  1]$, and then generate $y_i$ based on
%
%e5.1 #&#
\begin{equation}
\label{eq5.1} \mathsf{ P}(y_i = k) = \frac{\xi_i^{k}}{k!}
e^{-\xi_i},\qquad \log(\xi_i) = \sum_{j=1}^{d_n}
a_j(U_i) x_{ij}.
\end{equation}
We set the sample size $n=200$, the number of significant covariates $s_{n2}$
to be the integer part of $\ln n$ and $a_j(\cdot)$'s in (\ref{eq5.1})
to be
\begin{eqnarray*}
a_1(U)& =& -U, \qquad a_2(U) = \sin(2 \pi U),\qquad
a_3(U) = 4(U-0.5)^2,
\\
 a_4(U)& =& c_1= 0.6, \qquad a_5(U) =
c_2= -0.7, \qquad a_j(U)=0 \qquad\mbox {for } j > 5.
\end{eqnarray*}

For dimensions $d_n = 50$, $d_n = 100$, $d_n=200$, and $d_n = 500$, we
apply both the adaptive group LASSO method and the group SCAD method to
the simulated sample to select the model, and estimate the unknown
functional or constant coefficients. For each case, we do $1000$
simulations, and compute the mean integrated squared error
(MISE) of the estimators of the unknown functional coefficients and the
mean squared error (MSE) of the estimators of the unknown constant
coefficients. We also calculate the ratios of correct, under-selected,
under-specified, over-selected, over-specified and other models. The
``under-selected models'' means the selected models ignoring the
significant covariates. The ``under-specified models'' means where the
functional coefficients are mis-specified as the constant coefficients.
The ``over-selected models'' means the selected models including the
insignificant covariates. The ``over-specified models'' means where the
constant coefficients are mis-specified as functional. The ``other
models'' means that there exist more than one incorrect situation as
listed above. The ``correct models'' need not only select the true model
but also correctly identify the modelling structure.

%t1 #&#
\begin{table}[b]
\tabcolsep=0pt
\caption{The ratios of model selection in 1000 simulations}
\label{tab1}
\begin{tabular*}{\textwidth}{@{\extracolsep{\fill}}lcccccc@{}}
\hline
$\bolds{d_n}$ & \textbf{Correct} & \textbf{Under-selected} &
\textbf{Under-specified}
& \textbf{Over-selected} & \textbf{Over-specified} &
\textbf{Others}
\\
\hline
\multicolumn{7}{c}{Adaptive group LASSO method}
\\
\phantom{0}$50$ & 0.967 & 0.001 & 0.002 & 0.005 & 0.024 & 0.001
\\
$100$ & 0.944 & 0.001 & 0.005 & 0.008 & 0.039 & 0.003
\\
$200$ & 0.915 & 0.006 & 0.014 & 0.012 & 0.045 & 0.008
\\
$500$ & 0.863 & 0.015 & 0.023 & 0.026 & 0.057 & 0.016
\\[6pt]
\multicolumn{7}{c}{Group SCAD method}
\\
%$d_n$ & {\scriptsize Correct} & {\scriptsize Under-selected} &
%{\scriptsize Under-specified}
%& {\scriptsize Over-selected} & {\scriptsize Over-specified} &
%{\scriptsize Others}
%\\\hline
\phantom{0}$50$ & 0.970 & 0.001 & 0.002 & 0.003 & 0.022 & 0.002
\\
$100$ & 0.948 & 0.002 & 0.004 & 0.006 & 0.038 & 0.002
\\
$200$ & 0.925 & 0.005 & 0.012 & 0.010 & 0.042 & 0.006
\\
$500$ & 0.878 & 0.013 & 0.020 & 0.022 & 0.051 & 0.016
\\
\hline
\end{tabular*}
\tabnotetext[]{}{The ratios of choosing the correct, under-selected, under-specified,
over-selected, over-specified and other models in 1000 simulations by
using either the adaptive group LASSO method or the group SCAD method.}
\end{table}

The simulation results are reported in Tables~\ref{tab1} and \ref
{tab2}. We
can see from Table~\ref{tab1} that both the adaptive group LASSO method
and the
group SCAD method work well for model selection and structure
specification, and the group SCAD method gives slightly better
performance. Table~\ref{tab2} shows that the estimators obtained by
either the adaptive group LASSO method or the group SCAD method are
doing very well, and their performance is comparable to that of the
oracle estimators.
\end{exa}

\begin{exa}\label{ex5.2}
As the varying coefficient models are a
special case of the generalised varying coefficient models, our method
is also applicable to the varying coefficient models. In this example,
we compare our method with the KLASSO method proposed in \citet{WanXia09} for varying coefficient models. We consider exactly the same
simulated example as that in
\citet{WanXia09}, that is the following three varying coefficient models:
\begin{longlist}[(III)]
\item[(I)] $y_i=2 \sin(2 \pi U_i ) x_{i1} + 4U_i(1 - U_i) x_{i2} +
\sigma\varepsilon_i$,
\item[(II)] $y_i=\exp(2U_i-1) x_{i1} + 8U_i(1 - U_i) x_{i2} + 2\cos
^2(2 \pi U_i) x_{i3}
+\sigma\varepsilon_i$,
\item[(III)] $y_i=4U_i x_{i1} + 2 \sin(2 \pi U_i )x_{i2} + x_{i3} +
\sigma\varepsilon_i$,
\end{longlist}
where $x_{i1}=1$ for any $i$, $(x_{i2},  \ldots,  x_{i7})^{\T}$ and
$\varepsilon_{i}$, $i=1,  \ldots,  n$, are independently generated from
a multivariate normal distribution with $\operatorname{cov} (x_{ij_1},  x_{ij_2})=\break
0.5^{\llvert  j_{1}-j_{2}\rrvert  }$ for any $2\leq j_{1},  j_{2}\leq7$ and the standard normal distribution $\mathsf{ N}(0,  1)$,
respectively, $U_i$, $i=1,  \ldots,  n$, are independently
generated from either uniform distribution $\mathsf{ U}[0,  1]$ or Beta
distribution $\mathsf{ B}(4,  1)$, $\sigma$ is set to be $1.5$. For each
model, we conduct $200$ simulations, and in each simulation, we apply
either our method or the KLASSO to do model selection and estimation
and then make the comparison. We measure the performance of model
selection by reporting the percentages of correct-, under- and
over-fitting. The obtained results are presented in Table~\ref{tab3}.
From Table~\ref{tab3}, we can see our method performs better than the
KLASSO in model selection.

%t2 #&#
\begin{table}[t]
\caption{The MISEs and MSEs of the estimators for the functional and
constant coefficients}
\label{tab2}
\begin{tabular*}{\textwidth}{@{\extracolsep{\fill}}lccccccccc@{}}
\hline
\multicolumn{1}{c}{}
&\multicolumn{3}{c}{\textbf{Adaptive group LASSO}}
&\multicolumn{3}{c}{\textbf{Group SCAD}}
&\multicolumn{3}{c@{}}{\textbf{Oracle estimators}}
\\[-6pt]
\multicolumn{1}{c}{}
&\multicolumn{3}{c}{\hrulefill}
&\multicolumn{3}{c}{\hrulefill}
&\multicolumn{3}{c@{}}{\hrulefill}\\
$\bolds{d_n}$ & $\bolds{\widehat{a}_1(\cdot)}$ & $\bolds{\widehat{a}_2(\cdot)}$
& $\bolds{\widehat
{a}_3(\cdot)}$
& $\bolds{\overline{a}_1(\cdot)}$ & $\bolds{\overline{a}_2(\cdot)}$ & $\bolds{\overline
{a}_3(\cdot)}$
& $\bolds{a_{1}^{uo}(\cdot)}$ & $\bolds{a_{2}^{uo}(\cdot)}$ & $\bolds{a_{3}^{uo}(\cdot)}$
\\
\hline
\phantom{0}$50$ & 0.026 & 0.037 & 0.038 & 0.024 & 0.033 & 0.036 & 0.019 & 0.023 & 0.025
\\
$100$ & 0.038 & 0.049 & 0.050 & 0.035 & 0.045 & 0.048 & 0.019 & 0.023 & 0.025
\\
$200$ & 0.058 & 0.069 & 0.072 & 0.052 & 0.063 & 0.066 & 0.019 & 0.023 & 0.025
\\
$500$ & 0.090 & 0.095 & 0.098 & 0.084 & 0.093 & 0.091 & 0.019 & 0.023 & 0.025
\\
\hline
$\bolds{d_n}$ & $\bolds{\widehat{c}_1}$ & $\bolds{\widehat{c}_2}$ &
& $\bolds{\overline{c}_1}$ & $\bolds{\overline{c}_2}$ &
& $\bolds{c_{1}^{uo}}$ & $\bolds{c_{2}^{uo}}$ &
\\
\hline
\phantom{0}$50 $ & 0.015 & 0.017 & & 0.014 & 0.017 & & 0.006 & 0.008 &
\\
$100$ & 0.020 & 0.022 & & 0.018 & 0.021 & & 0.006 & 0.008 &
\\
$200$ & 0.030 & 0.035 & & 0.025 & 0.029 & & 0.006 & 0.008 &
\\
$500$ & 0.046 & 0.050 & & 0.039 & 0.042 & & 0.006 & 0.008 &
\\
\hline
\end{tabular*}
\tabnotetext[]{}{The MISEs or MSEs of the estimators obtained by either the adaptive
group LASSO method or the group SCAD method. For $j=1,2,3$ and $k=1,2$,
$\widehat{a}_j(\cdot)$'s and $\widehat{c}_k$'s are the
estimators obtained by the adaptive group LASSO method,
$\overline{a}_j(\cdot)$'s and $\overline{c}_k$'s are the
estimators obtained by the group SCAD method,
and $a_{j}^{uo}(\cdot)$'s and $c_{k}^{uo}$'s are the unbiased oracle estimators.}
\end{table}

%t3 #&#
\begin{table}
\caption{Comparison of model selection between our method and KLASSO}
\label{tab3}
\begin{tabular*}{\textwidth}{@{\extracolsep{\fill}}lccccccc@{}}
\hline
\multicolumn{2}{c}{}
&\multicolumn{3}{c}{\textbf{Our method}}
&\multicolumn{3}{c}{\textbf{KLASSO}}
\\[-6pt]
\multicolumn{2}{c}{}
&\multicolumn{3}{c}{\hrulefill}
&\multicolumn{3}{c}{\hrulefill}
\\
$\bolds{f_U(\cdot)}$ & $\bolds{n}$ & \textbf{Under} & \textbf{Correct} &
\textbf{Over}
& \textbf{Under} & \textbf{Correct} & \textbf{Over}
\\
\hline
\multicolumn{8}{c}{Model I}
\\
$\mathsf{ U}[0,1]$ & 100 & 0.020 & 0.910 & 0.070 & 0.09 & 0.74 & 0.16
\\
& 200 & 0.005 & 0.985 & 0.010 & 0.02 & 0.95 & 0.03
\\
$\mathsf{ B}[4,1]$ & 100 & 0.020 & 0.875 & 0.105 & 0.21 & 0.58 & 0.21
\\
& 200 & 0.005 & 0.950 & 0.045 & 0.08 & 0.86 & 0.05
\\[3pt]
\multicolumn{8}{c}{Model II}
\\
$\mathsf{ U}[0,1]$ & 100 & 0.015 & 0.915 & 0.070 & 0.01 & 0.83 & 0.16
\\
& 200 & 0.005 & 0.990 & 0.005 & 0.00 & 0.99 & 0.01
\\
$\mathsf{ B}[4,1]$ & 100 & 0.015 & 0.890 & 0.095 & 0.01 & 0.82 & 0.18
\\
& 200 & 0.005 & 0.970 & 0.025 & 0.00 & 0.96 & 0.04
\\[3pt]
\multicolumn{8}{c}{Model III}
\\
$\mathsf{ U}[0,1]$ & 100 & 0.010 & 0.935 & 0.055 & 0.02 & 0.85 & 0.13
\\
& 200 & 0.000 & 0.995 & 0.005 & 0.00 & 0.99 & 0.01
\\
$\mathsf{ B}[4,1]$ & 100 & 0.015 & 0.895 & 0.090 & 0.02 & 0.79 & 0.19
\\
& 200 & 0.005 & 0.975 & 0.020 & 0.00 & 0.96 & 0.04
\\
\hline
\end{tabular*}
\tabnotetext[]{}{The columns corresponding to ``Under'', ``Correct'' and ``Over'' are the
ratios of under-fitting, correct-fitting and over-fitting for our
method and KLASSO under different situations.}
\end{table}

As in \citet{WanXia09}, we employ the median of the relative
estimation errors (MREE), obtained in the $200$ simulations, to assess
the accuracy of an estimation method. The relative estimation error
(REE) is defined as
%
%e5.2 #&#
\begin{equation}
\label{eq5.2} \mathrm{REE} = 100 \times \frac{
\sum_{i=1}^{n} \sum_ {j=1}^{d_{n}}
\llvert
\hat{a}_{j}(U_{i}) - a_{j}(U_{i})
\rrvert }{
\sum_{i=1}^{n} \sum_{j=1}^{d_{n}}
\llvert
\hat{a}_{j}^{uo}(U_{i}) - a_{j}(U_{i}) \rrvert  },
\end{equation}
where $\hat{a}_{j}(\cdot)$ is the estimator of $a_j(\cdot)$, obtained
by the
estimation method concerned, and $\hat{a}_{j}^{uo}(\cdot)$ is the oracle
estimator of $a_j(\cdot)$.\vspace*{1pt} The median of REEs of our method and the KLASSO
under different situations
are presented in Table~\ref{tab4}, which shows our method is more
accurate than the KLASSO on estimation side. We thus conclude that our
method performs better than the KLASSO on both model selection and estimation.
%
%t4 #&#
\begin{table}
\caption{Comparison of estimation between our method and KLASSO}
\label{tab4}
\begin{tabular*}{\textwidth}{@{\extracolsep{\fill}}lccccccc@{}}
\hline
\multicolumn{8}{c}{\textbf{Median of relative estimation errors}}
\\
\hline
$\bolds{f_U(\cdot)}$ & $\bolds{n}$ & \textbf{Our method} &
\textbf{KLASSO} &$\bolds{f_U(\cdot)}$ & $\bolds{n}$ & \textbf{Our
method} & \textbf{KLASSO}
\\
\hline
\multicolumn{8}{c}{Model I}
\\
$\mathsf{ U}[0,1] $& 100 & 109.35 & 121.00 & $\mathsf{ B}[4,1]$ & 100 & 114.41 & 127.42
\\
$\mathsf{ U}[0,1] $& 200 & 101.78 & 115.45 & $\mathsf{ B}[4,1]$ & 200 & 103.49 & 122.12
\\[3pt]
\multicolumn{8}{c}{Model II}
\\
$\mathsf{ U}[0,1]$ & 100 & 107.81 & 109.45 & $\mathsf{ B}[4,1]$ & 100 & 115.17
& 111.06
\\
$\mathsf{ U}[0,1]$ & 200 & 101.51 & 109.46 & $\mathsf{ B}[4,1]$ & 200 & 103.73
& 108.07
\\[3pt]
\multicolumn{8}{c}{Model III}
\\
$\mathsf{ U}[0,1]$ & 100 & 106.71 & 116.53 & $\mathsf{ B}[4,1]$ & 100 & 112.39
& 118.91
\\
$\mathsf{ U}[0,1]$ & 200 & 101.21 & 110.59 & $\mathsf{ B}[4,1]$ & 200 & 104.06
& 113.43
\\
\hline
\end{tabular*}
\end{table}
%%
%%\begin{singlespace}
%%
%\begin{center}
%%
%\begin{minipage}{30pc}
%%{\it
%%The MREEs of our method or the KLASSO under different situations. The
%%column corresponding to $f_U(\cdot)$ is the distribution of $U$,
%%corresponding
%%to $n$ is the sample size. The columns corresponding to ``Our Method''
%%and
%%``KLASSO'' are the MREEs of our method and KLASSO, respectively.
%%}
%\end{minipage}
%%
%\end{center}
%%
%%\end{singlespace}
%%
%\end{table}
\end{exa}

\begin{exa}\label{ex5.3}
In this example, we compare the model
selection performance of our method with the method proposed in \citet{Lia12} for generalised varying coefficient models. We consider exactly
the same simulation settings as that in Example~2 of \citet{Lia12}, that
is the following varying coefficient logistic regression model where
the conditional mean function is
%
%e5.3 #&#
\begin{equation}
\label{eq5.3} \mathsf{ E}[y_i | X_i]=
\frac{\exp \{ \sum_{j=1}^{d_n} a_j(U_i)
x_{ij}  \}}{1+\exp \{ \sum_{j=1}^{d_n} a_j(U_i) x_{ij}
 \}}.
\end{equation}
The covariates are generated as following: for any $i=1,  \ldots,  n$, $x_{i1}=1$ and $(x_{i2},  \ldots,  x_{id_n})^{\T}$ are generated
from a
multivariate normal distribution with\break $\operatorname{cov} (x_{ij_1},  x_{ij_2})=0.1^{\llvert  j_{1}-j_{2}\rrvert  }$
for any $2 \leq j_{1},  j_{2} \leq d_n$. The index variable $U_i$,
$i=1,  \ldots,  n$, are independently
generated from the uniform distribution $\mathsf{ U}[0,  1]$. We set the
$a_j(\cdot)$'s in (\ref{eq5.3}) to be
\begin{eqnarray*}
a_1(U)& =& -4\bigl(U^3+2U^2-2U\bigr), \qquad
a_2(U) = 4\cos(2 \pi U), \\
a_3(U) &=& 3\exp\{U - 0.5 \}, \qquad
a_j(U) = 0 \qquad\mbox{when } j > 3.
\end{eqnarray*}

Similar to Example~2 of \citet{Lia12}, we set the sample size $n=150$ and
dimension $d_n=50$ or $d_n=200$. For each case, the simulation results
are based on 100 replications. The model selection performance is
measured by the average number of correct and incorrect varying
coefficients. The former one means the average number of significant
covariates that are correctly selected into the final model while the
latter means the average number of insignificant covariates that are
falsely selected as significant. The comparison results are shown in
Table~\ref{tab5}, from which we can see our method gives better model
selection results.

%t5 #&#
\begin{table}
\caption{Comparison between our method and Lian's methods on model selection}
\label{tab5}
\begin{tabular*}{\textwidth}{@{\extracolsep{4in minus 4in}}lcc@{}}
\hline
& \multicolumn{2}{c@{}}{\textbf{Average \# of varying coef.}}
\\[-6pt]
& \multicolumn{2}{c@{}}{\hrulefill}
\\
\textbf{Method} & \textbf{Correct} & \multicolumn{1}{c@{}}{\textbf{Incorrect}}
\\
\hline
& \multicolumn{2}{c}{$d_n=50$}
\\
GL(BIC) & 3\phantom{00.} & 18.75
\\
GL(eBIC) & 3\phantom{00.} & 16.33
\\
AGL(BIC-BIC) & 3\phantom{00.} & 10.29
\\
AGL(eBIC-eBIC)& 3\phantom{00.} & \phantom{0}1.56
\\
Our method & 3\phantom{00.} & \phantom{0}1.37
\\[3pt]
& \multicolumn{2}{c}{$d_n=200$}
\\
GL(BIC) & 3\phantom{00.} & 38.78
\\
GL(eBIC) & 3\phantom{00.} & 21.04
\\
AGL(BIC-BIC) & 3\phantom{00.} & 25.72
\\
AGL(eBIC-eBIC)& 2.96 & \phantom{0}2.49
\\
Our method & 3\phantom{00.} & \phantom{0}2.18
\\
\hline
\end{tabular*}
\tabnotetext[]{}{The simulation results are based on 100 replications with sample size
$n=150$. GL means group lasso method, AGL means adaptive group lasso
method. The details of GL and AGL methods can be found in \citet{Lia12}
and eBIC means extended Bayesian information criterion [\citet{CheChe08}].}
\end{table}
\end{exa}

%%%%%%%%%%%%%%%

%s6 #&#
\section{Real data analysis}
\label{sec6}

We now apply the proposed method to analyse an environmental data set
from Hong Kong. This data set was collected between January~1, 1994,
and December 31, 1995. It is a collection of numbers of daily total
hospital admissions for circulationary and respirationary problems,
measurements of pollutants and other environmental factors in Hong
Kong. The collected environmental factors are $\SO_2$ (coded by $x_1$),
$\NO_2$ (coded by $x_2$), dust (coded by $x_3$), temperature (coded by
$x_4$), change of temperature (coded by $x_5$),
humidity (coded by $x_6$) and ozone (coded by $x_7$). What we are
interested in is which environmental factors among the collected
factors have significant effects on the number of daily total hospital
admissions for circulationary and respirationary problems (coded by~$y$), and whether the impacts of those factors vary over time (coded by
$U$). As the numbers of daily total hospital admissions are count data,
it is natural to use Poisson regression model with varying
coefficients, namely (\ref{eq5.1}), to fit the data.

We apply the proposed group SCAD method to identify the significant
variables and the nonzero constant coefficients, and estimate the
functional or
constant coefficients in the selected model. The kernel function used
is still
taken to be the Epanechnikov kernel, and the bandwidth is chosen to be
$0.75[(\log d_n) / n]^{0.2}100 \%$ of the range of the time. The tuning
parameters $\lambda_{1}$, $\lambda_{2}$, $\lambda_{4}$ and $\lambda
_{4}^{\ast}$ are selected by the data driven approach described in
Section~\ref{sec4.2}.

\begin{figure}

\includegraphics{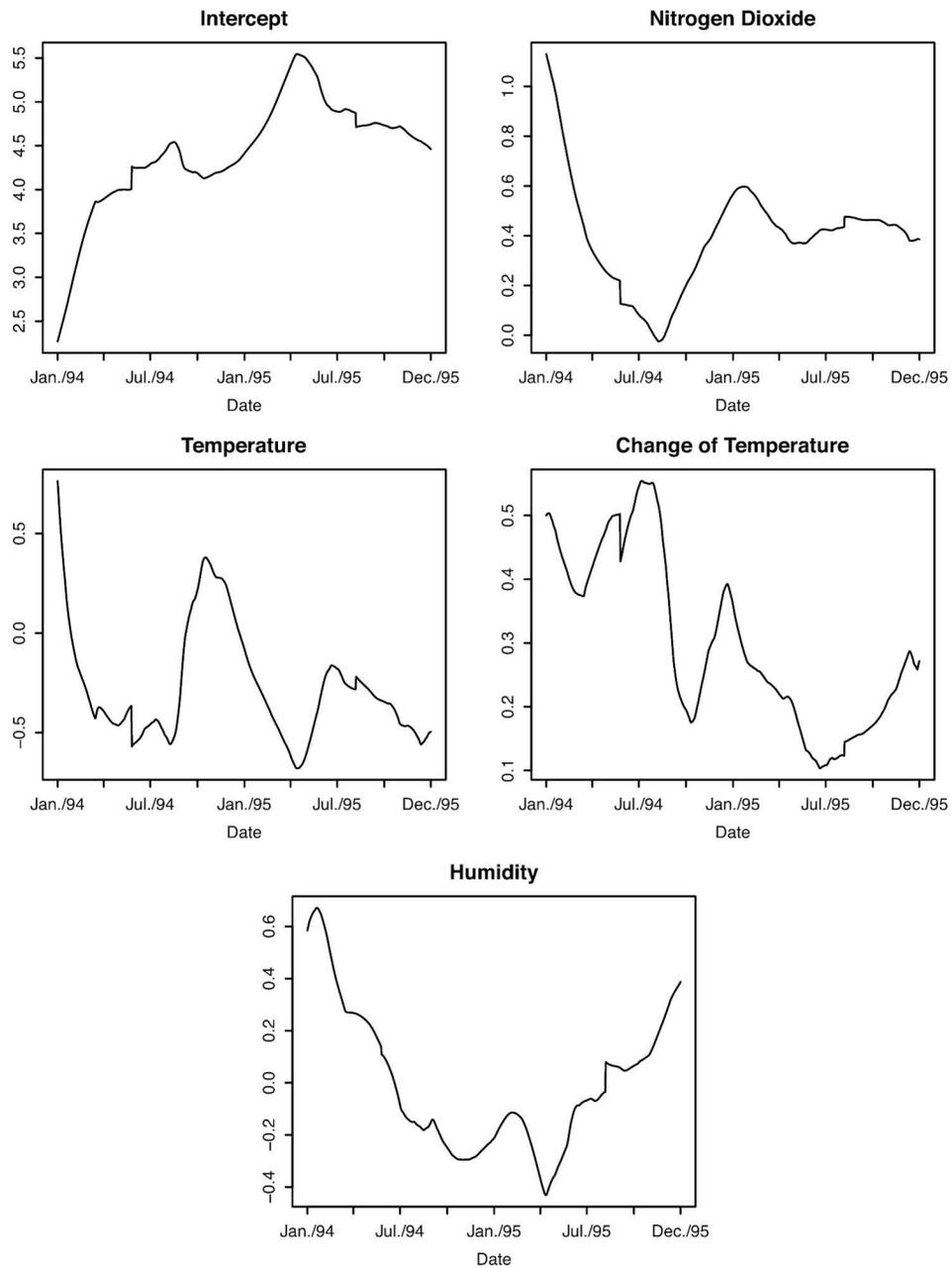}

\caption[Figure~1]{Estimated curves of the functional coefficients in
the selected model for the Hong Kong environment data.}
\label{fig1}
\end{figure}

The selected model is
\[
\mathsf{ P}(y_i = k) = \frac{\xi_i^{k}}{k!} e^{-\xi_i}
\]
with
\[
\log(\xi_i) = a_0(U_i) +
a_2(U_i) x_{i2} + a_4(U_i)
x_{i4} + a_5(U_i) x_{i5} +
a_6(U_i) x_{i6}.
\]
This shows only variables $\NO_2$, temperature, change of temperature
and humidity have effects on the number of daily total hospital
admissions for circulationary and respirationary problems, and all of
these variables have time-varying impacts. The estimates of the impacts
of these variables are presented in Figure~\ref{fig1}.

Figure~\ref{fig1} shows that $\mathrm{NO_{2}}$ always has a positive
impact on the daily number of total hospital admissions for
circulationary and respirationary problems, and this impact is stronger
in winter and spring than that in summer and autumn. This is in line
with the finding in one World Health Organization report [WHO report,
(\citeyear{WorOrg03})] which shows some evidence that ``long-term exposure
to $\mathrm{NO_{2}}$ at concentrations above 40--100 $\upmu\mathrm
{g}$/$m^{3}$ may decrease lung function and
increase the risk of respiratory symptoms.'' The nonlinear dynamic
pattern of the impact of $\mathrm{NO_{2}}$ also makes sense. This is
because the main source of $\mathrm{NO_{2}}$ pollution comes from the
burning of coals and gasolines. In the winter and spring season,
heating requirements will increase the amount of $\mathrm{NO_{2}}$
pollution. This is evident from the plot of $\mathrm{NO_{2}}$ in the
data set. Furthermore, the fog and mist in winter and spring will also
increase the chance that people expose to $\mathrm{NO_{2}}$. Though
$\mathrm{NO_{2}}$ is toxic by inhalation, as its compound is acrid and
easily detectable by smell at low concentrations, in most cases, the
inhalation exposure to $\mathrm{NO_{2}}$ can be generally avoided.
However, when $\mathrm{NO_{2}}$ is dissolved into the fog, this acid
mist will be hard to detect, and people may easily expose to this toxic
acid mist for a long time without awareness.

Figure~\ref{fig1} also shows the change of temperature has a
time-varying positive impact on the daily number of total hospital
admissions for
circulationary and respirationary problems. This coincides with the
intuition that a sudden change of temperature would greatly increase
the risk of catching cold, fever and other upper respirationary
diseases. The impact of temperature is also time varying and mostly
negative. It is stronger in autumn and spring than that in other
seasons. This makes sense, indeed, colder autumn or spring would see
more people catching circulationary or respirationary diseases.

%f1 #&#

The impact of humidity on the daily number of total hospital admissions
for circulationary and respirationary problems is interesting and
complicated. It does not seem to have any seasonal pattern. This is in
line with the findings reported in the literature. Indeed, existing
research [\citet{StrSan89,Sch95}; and \citet{deLetal96}] agrees that humidity has a significant effect on daily
hospital admissions for circulationary and respirationary problems in
many different places. \citet{StrSan89} study the childhood
respiratory problems against the indoor air temperature and relative
humidity. Through a randomly sampled questionnaire survey, and
interview of $1000$ children aged 7 about their living conditions and
reported circulationary and respirationary problems, they show that the
children living in damp (higher relative humidity level) bedrooms had
significantly higher probability to catch day cough, night cough and
chesty colds. \citet{Sch95} studies the short term fluctuations in
air pollution and hospital admissions of the elderly for
respiratory disease. According to their data set, the risk, measured by
sample variance, of respiratory hospital admissions of people aged 65
or above is bigger in the cities with higher average humidity levels
(measured by dew point). \citet{deLetal96} study the effects of
air pollution on daily hospital admissions for respiratory disease
based on a data set collected in London between 1987--1988 and
1991--1992. They show that the relative humidity is more significant for
the respiratory hospital admission numbers of children (0--14 years) and
the elderly (65$+$ years). All of these suggest that there may be a
strong relationship between humidity level and the risk for children
and elderly people to catch circulationary or respirationary disease.

Furthermore, we would like to examine the prediction performance of the
selected model and compare it with the full model with functional
coefficients. Given either model, we begin with using the first 700
observations as the training set to estimate the conditional
expectation of the response variable of the 701st observation. Then we
repeat this one-step forward prediction by enrolling one more
observation into the training set at a time. Finally, we end with using
the first 729 observations to predict the 730th observation. The
prediction performance is measured by the mean relative prediction
error (MRPE) defined as follows:
\[
\mathrm{MRPE}=\frac{1}{30}\sum_{i=701}^{730}
\biggl\llvert \frac{\hat{\xi
}_i-y_i}{y_i} \biggr\rrvert \times100\%,
\]
where $\hat{\xi}_i$ is the estimator of the conditional expectation of
the response variable at time $U_i$, $i=701, \ldots,730$. The MRPE
of the model selected by our method is $18.7\%$ while the MRPE of the
full model with functional coefficients is $41.6\%$. Hence, we can see
that the model selected by our method do have a better prediction
accuracy than the full model.

%\begin{figure}[htbp]
%\centerline{\psfig{figure=figure/Fig2_1.ps,width=5.0in}}
%\centerline{\psfig{figure=figure/Fig2_2.ps,width=5.0in}}
%\centerline{\psfig{figure=figure/Fig2_3.ps,width=5.0in}}
%\caption[Figure~1]{\it Estimates of the coefficients in the selected
%model for
%the Hong Kong environment data.}
%\label{fig1}
%\end{figure}

%\begin{figure}[htbp]
%\center{\includegraphics[width=5.0in]{Fig2_1.ps}}
%\center{\includegraphics[width=5.0in]{Fig2_2.ps}}
%\center{\includegraphics[width=5.0in]{Fig2_3.ps}}
%\caption[Figure~1]{Estimated curves of the functional coefficients in
%the selected model for the Hong Kong environment data.}
%\label{fig1}
%\end{figure}

%%%%%%%%%%%%%%%

%%%%%%%%%%%%%%%%%%%%%%

\begin{appendix}%\label{app}
\section{Assumptions}\label{app.A}
In this appendix, we give some regularity conditions which are
needed to prove the asymptotic theory. In Appendices B and C of the
supplemental material [Li, Ke and Zhang (\citeyear{supp})], we provide the proofs
of the main theoretical results and some auxiliary results, respectively.

%%%%%%%%%%%%%%%%%%%%%%%%%%%%%%%

%%s8 #&#
%
%\setcounter{equation}{0}
Recall that
\[
q_1(s,y)=\frac{\partial\ell [g^{-1}(s),y ]}{\partial s},\qquad q_2(s,y)=\frac{\partial^2 \ell [g^{-1}(s),y ]}{\partial s^2}
\]
and define
\[
\ddot{\mathcal{L}}_{n}(u)=\left[ %
\matrix{
\ddot{\mathcal{L}}_{n}(u, 0)&\ddot{\mathcal{L}}_{n}(u,1)
\vspace*{2pt}\cr
\ddot{\mathcal{L}}_{n}(u, 1)&\ddot{\mathcal{L}}_{n}(u, 2)}
 \right]
\]
with
\[
\ddot{\mathcal{L}}_{n}(u,l) =\frac{1}{n}\sum
_{i=1}^n q_2 \Biggl[\sum
_{j=1}^{d_n}a_j(U_i)x_{ij},
y_i \Biggr] X_iX_i^{\T} \biggl(
\frac{U_i-u}{h} \biggr)^{l}K_h(U_i - u)
\]
for $l=0,1,2$. Define $b= \max\{\lambda_1/\lambda_2, \lambda
_2/\lambda
_1\}+\delta$ for any $\delta>0$, where $\lambda_1$ and $\lambda_2$ are
defined in (\ref{eq2.3}), and let
\begin{eqnarray*}
{\bolds\Omega}_0(b)&=& \Biggl\{  \mathbf{v}=(v_{11},
\ldots ,v_{1d_n},v_{21},\ldots,v_{2d_n})^{\T}:
\Vert{\mathbf v}\Vert=1,
\\
&& \sum_{j=1}^{d_n}\bigl(|v_{1j}|+|v_{2j}|\bigr)
\leq2(1+b) \sum_{j=1}^{s_{n2}}\bigl(|v_{1j}|+|v_{2j}|\bigr)
\Biggr\}.
\end{eqnarray*}
When $\lambda_{1}\propto\lambda_{2}$ (see Assumptions \ref{ass5} or
\ref{ass5prime}
below), $b$ is bounded by a positive constant, and it becomes $1+\delta
$ which could be sufficiently close to $1$ if we further assume that
$\lambda_1=\lambda_2$. To simplify the presentation, we denote
\[
Q_{i1}=q_1 \Biggl[\sum_{j=1}^{d_n}a_j(U_{i})x_{ij},
y_i \Biggr],\qquad Q_{i2}=q_2 \Biggl[\sum
_{j=1}^{d_n}a_j(U_{i})x_{ij},
y_i \Biggr].
\]
We next introduce some regularity conditions which are needed to
establish the asymptotic theory for the proposed model selection and
structure specification procedure. Some of the conditions might be not
the weakest possible conditions.

%%
%\begin{description}
\begin{ass}\label{ass1}
The kernel function $K(\cdot)$ is a
continuous and symmetric probability density function with a compact support.
\end{ass}

\begin{ass}\label{ass2} (i) Let $\mathsf{ E}
(Q_{i1} |
X_i, U_i )=0$ a.s.,
and ${\mathsf{E}} (Q_{i1}^2 |U_i=u )$ be continuous for
$u\in
[0, 1]$. Moreover, suppose that uniformly for $u\in[0,1]$, either
%
%e8.1 #&#
\begin{eqnarray}
\label{eqA.1}&&\max_{1\leq j\leq d_n}\mathsf{ E} \bigl[|Q_{i1}x_{ij}|^{m_0}
|U_i=u \bigr]+\max_{1\leq j,k\leq d_n}\mathsf{ E}
\bigl[|Q_{i2}x_{ij}x_{ik}|^{m_0}
|U_i=u \bigr]
\nonumber
\\[-8pt]
\\[-8pt]
\nonumber
&&\qquad<\infty\qquad \mbox{a.s.}
\end{eqnarray}
for $m_0>2$, or
%
%e8.2 #&#
\begin{eqnarray}
\label{eqA.2}&& \max_{1\leq j\leq d_n}\mathsf{ E} \bigl[|Q_{i1}x_{ij}|^{m}
| U_i=u \bigr]+\max_{1\leq j,k\leq d_n}\mathsf{ E}
\bigl[|Q_{i2}x_{ij}x_{ik}|^{m}
|U_i=u \bigr]
\nonumber
\\[-8pt]
\\[-8pt]
\nonumber
&&\qquad\leq\frac{M_0 m!}{2} \qquad\mbox{a.s.}
\end{eqnarray}
for all $m\geq2$ and $0<M_0<\infty$.

(ii) Let $q_2(s,y)<0$ for $s\in{\mathbb{R}}$ and $y$ in the
range of the response variable. Furthermore, there exists an
$M(X,U,y)>0$ such that
\[
\bigl|q_2 \bigl[r(X,U)+\delta_\ast,y \bigr]-q_2
\bigl[r(X,U),y \bigr]\bigr |\leq M(X,U,y)|\delta_\ast|
\]
with $r(X,U)=\sum_{j=1}^{d_n}a_j(U)x_{j}$, and uniformly for $u\in
[0,1]$ either
\[
\max_{1\leq j,k,l\leq d_n}{\mathsf{E}}
\bigl[\bigl|x_{ij}x_{ik}x_{il}M(X_i,U_i,y_i)\bigr|^{m_0}
| U_i=u \bigr]<\infty \qquad \mbox{a.s.}
\]
for $m_0>2$ if (\ref{eqA.1}) is satisfied, or
\[
\max_{1\leq j,k,l\leq d_n}{\mathsf{E}}
\bigl[\bigl|x_{ij}x_{ik}x_{il}M(X_i,U_i,y_i)\bigr|^{m}|U_i=u
\bigr]<\frac{M_1 m!}{2}\qquad \mbox{a.s.}
\]
for all $m\geq2$ if (\ref{eqA.2}) is satisfied, $0<M_1<\infty$.

(iii) There exist $0<\rho_1\leq\rho_2<\infty$ such that
\[
\rho_1\leq\inf_{u\in[0,1]}\inf_{{\mathbf v}\in{\bolds\Omega
}_0(b)}{
\mathbf v}^{\T} \bigl[-\ddot{\mathcal{L}}_{n}(u) \bigr]{
\mathbf v}\leq\sup_{u\in
[0,1]}\sup_{{\mathbf v}\in{\bolds\Omega}_0(b)}{\mathbf
v}^{\T} \bigl[-\ddot {\mathcal{L}}_{n}(u) \bigr]{\mathbf v}
\leq\rho_2
\]
with probability approaching one.
\end{ass}

\begin{ass}\label{ass3}
The density function $f_U(\cdot)$ has a
continuous second-order derivative. In addition, $f_U(u)$ is bounded
away from zero and infinity when $u\in[0, 1]$.
\end{ass}

\begin{ass}\label{ass4}The functional coefficients,
$a_j(\cdot)$,
have continuous second-order derivatives for $j=1,\ldots,d_n$.
\end{ass}

\begin{ass}\label{ass5}Let $d_n\propto n^{\tau_1}$ and
$\frac
{nh}{(nd_n^3)^{2/m_0}\log h^{-1}}\rightarrow\infty$, where $0\leq
\tau
_1<\infty$ and $m_0$ is defined in (\ref{eqA.1}). Moreover, the
bandwidth $h$ and the tuning parameters $\lambda_1$ and $\lambda_2$
satisfy $h\propto n^{-\delta_1}$ with $0<\delta_1<1$, $s_{n2}h^2+
(\frac{\log h^{-1}}{nh} )^{1/2}=o(\lambda_{1})$, $\lambda
_{1}\propto
\lambda_{2}$ and $s_{n2}\lambda_1^2h^{-2}+s_{n2}^2\lambda_1=o(1)$.
\end{ass}

{\renewcommand{\theass}{5$^\prime$}
\begin{ass}\label{ass5prime}
%\item\textit{Assumption $\it5 }.
Let $d_n\propto\exp
 \{
(nh)^{\tau_2} \}$ with $0\leq\tau_2<1-\tau_3$, $0<\tau_3<1$.
Furthermore, the bandwidth $h$ and the tuning parameters $\lambda_1$
and $\lambda_2$ satisfy $h\propto n^{-\delta_1}$ with $0<\delta_1<1$,
$s_{n2}h^2+ (\frac{\log h^{-1}}{nh} )^{\tau_3/2}=o(\lambda_{1})$,
$\lambda_{1}\propto\lambda_{2}$ and $s_{n2}\lambda
_1^2h^{-2}+s_{n2}^2\lambda_1=o(1)$.
\end{ass}}

\setcounter{ass}{5}
\begin{ass}\label{ass6}
(i) Let $s_{n2}h^2\propto
(nh)^{-1/2}$, $\lambda_4\sim\lambda_4^\ast$, $\lambda_4=o
(n^{1/2} )$, and
%
%e8.3 #&#
\begin{equation}
\label{eqA.3} h^{-1/2} \bigl[\bigl(\log h^{-1}
\bigr)^{1/2}+s_{n2}^{1/2}(1+\lambda_1\sqrt
{nh}) \bigr]=o(\lambda_4).
\end{equation}

(ii) There exists a positive constant $b_\diamond$ such that
%
%e8.4 #&#
\begin{equation}
\label{eqA.4} \min_{1\leq j\leq s_{n2}}\|{\balpha}_{j0}\|\geq
b_\diamond n^{1/2},\qquad \min_{1\leq j\leq s_{n1}}D_j
\geq b_\diamond n^{1/2}
\end{equation}
with probability approaching one. Furthermore, uniformly for
$k=s_{n2}+1,\ldots,\break d_n, d_n+s_{n1}+1,\ldots,2d_n$,
%
%e8.5 #&#
\begin{equation}
\label{eqA.5} \sup_{u\in[0,1]} \sup_{\|{\mathbf w}^o\|=1} \bigl
\llvert \ddot{\mathcal{L}}_{n}(u | k) {\mathbf w}^o\bigr
\rrvert =O_P(1),
\end{equation}
where $\ddot{\mathcal{L}}_{n}(u | k)$ is the $k$th row of $\ddot
{\mathcal{L}}_{n}(u)$, ${\mathbf w}^o= [({\mathbf w}_1^o)^{\T}, ({\mathbf
w}_2^o)^{\T} ]^{\T}$, ${\mathbf w}_1^o$ and ${\mathbf w}_2^o$ are two
$d_n$-dimensional column vectors, the last $d_{n}-s_{n2}$ elements of
${\mathbf w}_1^o$ and the last $d_{n}-s_{n1}$ elements of ${\mathbf w}_2^o$ are zeros.
\end{ass}

\begin{rem}\label{rea1}
The above assumptions are mild and
justifiable. Assumption~\ref{ass1} is a commonly-used condition on the kernel
function and can be satisfied for the uniform kernel function and the
Epanechnikov kernel function which is used in our numerical study. The
compact support restriction on the kernel function is not essential and
can be removed at the cost of more tedious proofs. Assumption~\ref{ass2} imposes
some smoothness and moment conditions on $Q_{i1}$ and $Q_{i2}$, some of
which are commonly used in local maximum likelihood estimation [cf.,
Cai, Fan and Li (\citeyear{CaiFanLi00}), \citet{LiLia08}]. Two moment conditions
(\ref
{eqA.1}) and (\ref{eqA.2}) are imposed in Assumption~\ref{ass2}(i), and they are
used to handle the polynomially diverging dimension of the covariates
(in Assumption~\ref{ass5}) and the exponentially diverging dimension of the
covariates (in Assumption~\ref{ass5prime}), respectively. Hence, as the
dimension of the covariates increase from the polynomial order to the
exponential order, the required moment condition would be stronger. In
contrast, most of the existing literature such as \citet{Lia12} only
considers the case of the stronger moment condition in (\ref{eqA.2}),
which may possibly limit the applicability of the model selection
methodology. Assumption~\ref{ass2}(iii) can be seen as the modified version of
the so-called \textit{restricted eigenvalue condition} introduced by
Bickel, Ritov and
Tsybakov (\citeyear{BicRitTsy09}) for the parametric regression models.
Assumptions \ref{ass3} and \ref{ass4} provide some smoothness conditions on the density
function of $U$ and the functional coefficients $a_j(\cdot)$, which are
not uncommon when the local linear approach is applied [cf., \citet{FanGij96}].

Assumption~\ref{ass5} imposes some restrictions on the bandwidth $h$ and the
tuning parameters $\lambda_1$ and $\lambda_2$ when $d_n\propto
n^{\tau
_1}$, whereas Assumption~\ref{ass5prime} imposes some conditions when
$d_n\propto\exp \{(nh)^{\tau_2} \}$. They are crucial to derive
the uniform convergence rates for the preliminary estimation in
Proposition~\ref{pr3.1}. Consider a special case when $s_{n2}$ is a fixed
positive integer; we may choose $h\propto n^{-1/5}$ and $\lambda
_1=\lambda_2\propto n^{-3/10}$. Then the conditions in Assumption~\ref{ass5}
would be satisfied when $m_0$ is sufficiently large, and those in
Assumption~\ref{ass5prime} would be satisfied when $3/4<\tau_3<1$. Noting
that $s_{n2}h^2+ (\frac{\log h^{-1}}{nh} )^{1/2}=o(\lambda_{1})$
and $\lambda_{1}\propto\lambda_{2}$ by Assumption~\ref{ass5}, the influence by
$h$ and $\lambda_2$ on the uniform convergence rate in (\ref{eq3.1}) is
dominated by that of $\lambda_1$. The regularity conditions in
Assumption~\ref{ass6} is mainly used to prove the sparsity and oracle property
for the proposed model selection procedure. By Assumption~\ref{ass5}, we may
show that the leading term of the left-hand side of (\ref{eqA.3}) is
$\lambda_1\sqrt{ns_{n2}}$. Once again we consider the special case when
$s_{n2}$ is a fixed positive integer and $h\propto n^{-1/5}$. Then the
conditions in Assumption~\ref{ass6} would be satisfied if we choose $\lambda
_1\propto n^{-3/10}$ and $\lambda_4=\lambda_4^\ast\propto n^{\tau_4}$
with $1/5<\tau_4<1/2$.
\end{rem}
\end{appendix}

%s7 #&#
\section*{Acknowledgements}
\label{sec7}
The authors thank the Co-Editor, an Associate Editor and three referees
for the helpful comments which greatly improved the former version of
the paper.

\begin{supplement}[id=suppA]
\stitle{Supplement to ``Model selection and structure specification in
ultra-high dimensional generalised semi-varying coefficient models''\\}
\slink[doi]{10.1214/15-AOS1356SUPP} %[doi,text={...}] - jei reikia suskaldyti doi
\sdatatype{.pdf}
\sfilename{aos1356\_supp.pdf}
\sdescription{We provide the detailed proofs of the main
results stated in Section~\ref{sec3}
as well as some technical lemmas which are useful in the proofs.}
\end{supplement}

%%%%%%%%%%%%%%%%%%%%%%%%%%%%%%%

% imsref loaded by akundreckaite, 2015-07-16 12:43:38
%

%\begin{appendix}
%\section{}
%\end{appendix}

%\section*{Acknowledgments}

%\begin{thebibliography}{99}
%\bibitem[\protect\citeauthoryear{}{}]{r1}
%\bibitem{r1}

\printaddresses
\end{document}